\documentclass[10pt, a4paper]{article}
\usepackage{graphicx}
\usepackage{amssymb, amsmath}
\usepackage{color}
\usepackage{grffile}
\usepackage{epstopdf}
\usepackage[ruled,vlined]{algorithm2e}
\usepackage[lofdepth,lotdepth, caption=false]{subfig}
\usepackage{dsfont}

\newcommand{\R}{\mathds{R}}
\newcommand{\id}{\mathds{1}}

\newcommand{\E}{\ensuremath{\mathbb{E}}}
\newcommand{\Tr}{\mathrm{Tr}}

\newcommand{\x}{\mathbf{x}}
\newcommand{\xx}{\mathbf{x}}
\newcommand{\yy}{\mathbf{y}}
\newcommand{\pp}{\mathbf{p}}

\newcommand{\uu}{\mathbf{u}}
\newcommand{\vv}{\mathbf{v}}

\DeclareMathOperator*{\argmin}{arg\,min}

\begin{document}
\pagestyle{headings}  

\title{On low complexity Acceleration Techniques for Randomized Optimization:\\ Supplementary Online Material}
%
%
\author{Sebastian U. Stich\footnotemark[1]}
%
%
%

\date{June 8, 2014}

\renewcommand{\thefootnote}{\fnsymbol{footnote}}
\footnotetext[1]{Institute of Theoretical Computer Science, ETH Z\"{u}rich, \texttt{sstich@inf.ethz.ch}}
\renewcommand{\thefootnote}{\arabic{footnote}}

\maketitle              
\begin{abstract}
Recently it was shown by Nesterov (2011) that techniques form convex optimization can be used to successfully accelerate simple derivative-free randomized optimization methods. The appeal of those schemes lies in their low complexity, which is only $\Theta(n)$ per iteration---compared to $\Theta(n^2)$ for algorithms storing second-order information or covariance matrices. From a high-level point of view, those accelerated schemes employ correlations between successive iterates---a concept looking similar to the evolution path used in Covariance Matrix Adaptation Evolution Strategies (CMA-ES). 
In this contribution, we (i) implement and empirically test a simple accelerated random search scheme (SARP). Our study is the first to provide numerical evidence that SARP can effectively be implemented with adaptive step size control and does not require access to gradient or advanced line search oracles. We (ii) try to empirically verify the supposed analogy between the evolution path and SARP. We propose an algorithm CMA-EP that uses only the evolution path to bias the search. 
This algorithm can be generalized to a family of low memory schemes, with complexity $\Theta(mn)$ per iteration, following a recent approach by Loshchilov (2014).  The study shows that the performance of CMA-EP heavily depends on the spectra of the objective function and thus it cannot accelerate as consistently as SARP.
~
\end{abstract}
\section{Introduction}
\label{sec:intro}
The Gradient Method~\cite{Polyak:1987,Nesterov:2004}---one of the most fundamental schemes in convex optimization---has iteration complexity $\Theta(n)$, where $n$ is the dimension. On strongly convex functions its convergence rate is linear, depending only on the condition number of the objective function. To overcome the difficulty imposed by ill-conditioned problems, second-order methods like Newton's method 
or first order Quasi-Newton methods such as the BFGS scheme \cite{Broyden:1970,Fletcher:1970,Goldfarb:1970,Shanno:1970} are a welcome alternative.
Those schemes maintain a quadratic model of the objective function and their complexity is bounded by $\Omega(n^2)$. Limited memory schemes like L-BFGS~\cite{Nocedal:1980,Liu:1989fu} trade-off linear iteration complexity $\Theta(mn)$ (where $m$ is a fixed parameter), versus convergence rate. Accelerated versions of the Gradient Method have linear complexity $\Theta(n)$ per iteration and converge with optimal rate among all first-order methods. On strongly convex problems the convergence rate is  proportional to the square root of the condition number~\cite{Polyak:1987,Nesterov:2004,nesterov:1983,nesterov:2007,tseng:2008}.

Randomized (gradient-free) schemes do not require first-order information, they operate by only querying function values. Such schemes are nowadays a ubiquitous tool for solving many practical problems in science and engineering  where first-order information is difficult to compute or does not exist. Among the first proposed schemes that are still of considerable (theoretical) importance are Adaptive Step Size Random Search (aSS)~\cite{Schumer:1968} and the (almost identical) well-known (1+1)-Evolution Strategy (ES)~\cite{Rechenberg:1973} in Evolutionary Computation. More recent schemes comprise Random Pursuit (RP)~\cite{Mutseniyeks:1964,Stich:2011}, or Random Gradient Descent~\cite{nesterov:11}. Those schemes can be viewed as generalizations of the Gradient Method to zeroth-order, with iteration complexity $\Theta(n)$. Likewise, analogues of the second-order schemes try to estimate an approximation of the Hessian by finite difference computations~\cite{Leventhal:2011,Stich:2012} or by estimating correlations among search directions. A very popular algorithm of this kind is the Covariance Matrix Adaptation Evolution Strategy (CMA-ES)~\cite{Hansen:2001,Hansen:2003}. Limited memory variants have been proposed in~\cite{knight:2007,Loshchilov:2014}, with iteration complexity $\Theta(mn)$. Especially, the later variant due to Loshchilov shows excellent convergence in high dimensions also for small values of $m$.
Only recently, also zeroth-order analogues of the accelerated gradient schemes have been introduced~\cite{nesterov:11,lee:2013}. Those schemes massively outperform the simple random search schemes on convex problems. This performance gain does not come for free, as those schemes require valid bounds on the condition number as input parameters. However, their low iteration complexity of $\Theta(n)$ could make them a promising choice for large scale problems, where the fully-quadratic schemes inherently fail. We focus on a very simple accelerated random search scheme, which we call SARP. 

By inspecting closely the accelerated search schemes, one could conclude that the difference to the classical schemes can be explained by an additional ``drift'' term~\cite{Nesterov:2004,tseng:2008} that takes into account correlations of the last iterates. In Evolutionary Computation, correlations between successive iterates are often expressed by the evolution path~\cite{ostermeier1994}. In the popular CMA-ES, the evolution path accumulates the correlation of successive iterates over a finite horizon of the order of $\Theta(n)$ steps~\cite{Hansen:2001}. In this work, we are interested, if the evolution path can be used for acceleration, competitive to the accelerated zeroth-order schemes form convex optimization. 
To this end, we introduce a variant of CMA-ES, called EP-CMA, that only uses the information stored in the evolution path to bias the direction of the search.
Similar to the approach proposed in~\cite{Loshchilov:2014}, this scheme can be generalized to a family of schemes, which we call EP-CMA-$m$. 
 For $m=1$, the approach is similar to~\cite{sun:2013}. For $m>1$, the schemes are similar in spirit to the LM-BFGS schemes~\cite{Nocedal:1980,Liu:1989fu} and for $m$ large enough, the scheme is intended to approach $(1+1)$-CMA-ES. EM-CMA-$m$ could be implemented with $\Theta(mn)$ complexity per iteration (not presented here), although this was not required, as the dimension $n \leq 100$ in our empirical study.

The remainder of this paper is structured as follows. In Section~\ref{sec:algorithms} we present the accelerated random search scheme SARP and detail the EP-CMA-$m$ schemes. In Section~\ref{sec:emp} we empirically test the performance of all schemes on three quadratic and the non-convex Rosenbrock function, and highlight the key results. We discuss these results and conclude the paper in Section~\ref{sec:disc}.
%
%
\def \sfactor{0.83}
 \begin{figure}[!ht]
    \centering
    \null\hfill %
     \scalebox{\sfactor}{
     \begin{minipage}[t]{.55\linewidth}
     \vspace{0pt}
     \begin{algorithm}[H]
     %
     %
      \NoCaptionOfAlgo
      \SetAlgoVlined
      \DontPrintSemicolon
      \SetInd{0.1em}{0.5em}
      \nl \lIf{exact}{
          \Return $\mathrm{exactLS}(\x, \uu/\|\uu\|)$
        }
        \nl \lElse{
           \Return $\mathrm{aSS}(\x, \uu, \sigma, p)$
        }%
      \caption{{lineSearch$(\x, \uu, [\sigma, p])$}}%
     \end{algorithm}%
     \begin{algorithm}[H]
     %
     %
      \NoCaptionOfAlgo
      \SetAlgoVlined
      \DontPrintSemicolon
      \SetInd{0.1em}{0.5em}
      \nl $\sigma_+ \leftarrow \min_{\lambda}f(\x+\lambda \uu);\, \x_+ \leftarrow \x + \sigma_+ \uu$\;
      \nl \Return{$(\x_+, \sigma_+)$}
      \caption{{exactLS$(\x, \uu, [\sigma])$}}%
     \end{algorithm}%
    \end{minipage}%
    }
    \hfill
    \scalebox{\sfactor}{
    \begin{minipage}[t]{.55\linewidth}
    \vspace{0pt}
        \begin{algorithm}[H]
     %
     %
      \NoCaptionOfAlgo
      \SetAlgoVlined
      \DontPrintSemicolon
      \SetInd{0.1em}{0.5em}
      \nl \eIf{ $f(\x+ \sigma \uu) \leq f(\x)$}{
      \nl $\x_+ \leftarrow \x + \sigma \uu$; $\sigma_+ \leftarrow \sigma \cdot \exp(1/3)$\;
      }{
      \nl $\x_+ \leftarrow \x$; $\sigma_+ \leftarrow \sigma \cdot \exp\left(-\frac{p}{3(1-p)}\right)$\;
      }%
      \nl \Return{$(\x_+, \sigma_+)$}\;%
      \caption{{aSS$(\x, \uu, \sigma,p)$ \textit{(adaptive step size)}}}%
     \end{algorithm}%
    \end{minipage}%
    }
     \hfill\null\par
    \caption{Line search oracles for gradient-free optimization.}\label{algo:ls}
   \end{figure} 
\section{Algorithms}
\label{sec:algorithms}
We here present the optimization schemes considered in this study. We first detail two Random Pursuit algorithms and a simplistic variant of a standard (1+1)-CMA-ES. Then we introduce the new EP-CMA-$m$ schemes.

\textbf{RP.}
Random Pursuit is a basic optimization scheme that iteratively generates a sequence of approximate solutions to the global optimization problem $\min_{\xx \in \R^n} f(\xx)$. In each step a search direction is drawn $\uu_k \sim \mathcal{N}(0, \id_n)$. In Random Pursuit with exact line search (RP-exact), first proposed in~\cite{Mutseniyeks:1964} and analyzed in~\cite{Stich:2011}, the step size $\sigma$ is determined by minimizing the objective function in direction $\uu$, i.e. $\sigma = \argmin_\lambda f(\x + \lambda \uu)$. For quadratic functions $f(\x):= \frac{1}{2} \x^T A \x$ with Hessian $A$, the expected one-step progress can be estimated as:
\begin{align}
\label{eq:generalRP}
\E \left[f(\x_+) \mid \x\right] \leq \left(1- \lambda_{\rm min}(A)/\Tr[A] \right) f(\x) \,,
\end{align}
where $\x$ is the current iterate, $\x_+ :=\x + \sigma \uu$ denotes the next iterate. This statement can also be generalized to arbitrary smooth convex functions \cite{Stich:2011}.
Stich et al. \cite{Stich:2011} show that RP-exact still converges if the line search is not performed exactly, but allowing relative errors. Therefore, we also consider Random Pursuit with adaptive step sizes (RP) instead of exact line search. In RP the step size is dynamically controlled such as to approximately guarantee a certain probability $p$ of finding an improving iterate. Depending on 
the underlying test function, different optimality conditions can be formulated for the value $p$. 
Schumer and Steiglitz~\cite{Schumer:1968} suggest the setting $p = 0.27$ which is considered throughout this work. We use immediate exponential step size control as explicitly formulated in the aSS sub-routine in Fig.~\ref{algo:ls}. RP is identical to the well known (1+1)-ES. 
%
%
\def \sfactor{0.83}
 \begin{figure}[!ht]
    \centering
    \null\hfill %
     \scalebox{\sfactor}{
     \begin{minipage}[t]{.55\linewidth}
     \vspace{0pt} 
     \begin{algorithm}[H]
     %
     %
      \NoCaptionOfAlgo
      \SetAlgoVlined
      \DontPrintSemicolon
      \SetInd{0.1em}{0.5em}
      \nl \For{$k=1$ to $N$}{
       \nl $\uu_{k} \sim \mathcal{N}(0, I_n)$\;
      \nl $(\x_{k},\sigma_k) \leftarrow \mathrm{lineSearch}(\x_{k-1}, \uu_{k}, [\sigma_{k-1}])$\;
      }%
      \nl \Return{$\x_N$}\;
      \caption{{RP$(\x_0, N, [\sigma_0 , p])$}}%
     \end{algorithm}%
      \begin{algorithm}[H]
     %
     %
      \NoCaptionOfAlgo
      \SetAlgoVlined
      \DontPrintSemicolon
      \SetInd{0.1em}{0.5em}
      \nl $\hat{\pp}_0 \leftarrow \mathbf{0};\, \hat{\pp}_{1},\dots,\hat{\pp}_{m-1} \leftarrow \mathbf{0};\,q=0$\;
      \nl \For{$k=1$ to $N$}{
       \nl $C_k \leftarrow I_n$\; 
       \nl \For{$i =1$ to $m-1$}{$C_{k} \leftarrow (1-c_{\rm{cov}}) C_{k}+c_{\rm{cov}} \hat{\pp}_{i} \hat{\pp}_{i}^T$}
       \nl $C_{k} \leftarrow (1-c_{\rm{cov}}) C_{k}+c_{\rm{cov}} \pp_{k-1} \pp_{k-1}^T$
       \nl $\uu_{k} \sim \mathcal{N}(0, C_{k})$\;     
       \nl $(\x_k, \sigma_k) \leftarrow \mathrm{aSS}(\x_{k-1}, \uu_k, \sigma_{k-1})$\;
       \nl $\yy_{k} \leftarrow (\xx_k - \xx_{k-1})/\sigma_{k-1}$\;
       \nl \uIf{$\yy_{k} \neq \mathbf{0}$ (\textit{success})}{
       \nl $\pp_{k} \leftarrow (1-c_{\rm{c}}) \pp_{k-1} + \sqrt{c_{\rm{c}}(2-c_{\rm{c}})}\yy_{k}$ \;
       }
       \nl \lElse{$\pp_{k} \leftarrow (1-c_{\rm p})\pp_{k-1}$}
       \nl \If{$k > q + n^2/m$}{
       $\hat{\pp}_{1} \leftarrow \hat{\pp}_{2}, \dots, \hat{\pp}_{m-2} \leftarrow \hat{\pp}_{m-1};\, q=k$
       }
      }%
      \nl \Return{$\x_N$}\;
      \caption{{EP-CMA-$m$ $(\x_0, N, \sigma_0, p, c_{\rm{c}}, c_{\rm{cov}})$}}%
     \end{algorithm}%
    \end{minipage}%
    }
    \hfill
    \scalebox{\sfactor}{
    \begin{minipage}[t]{.55\linewidth}
    \vspace{0pt} 
      \begin{algorithm}[H]
     %
     %
      \NoCaptionOfAlgo
      \SetAlgoVlined
      \DontPrintSemicolon
      \SetInd{0.1em}{0.5em}
      \nl $\yy_0 \leftarrow \x_0;\, \vv_0 \leftarrow x_0;\, \theta \leftarrow \sqrt{\frac{m}{2 n^2 L}}$\;
      \nl \For{$k=1$ to $N$}{
       \nl $\uu_{k} \sim \mathcal{N}(0, I_n)$\;
      \nl $(\x_{k},\sigma_k) \leftarrow \mathrm{lineSearch}(\yy_{k-1}, \uu_{k}, [\sigma_{k-1}])$\;
       \nl $\yy_k \leftarrow (\theta \vv_{k-1} + \x_{k}) /(1+\theta)$\;
       \nl $\vv_k \leftarrow (1-\theta)\vv_{k-1} + \theta \yy_{k} + \theta n \frac{L}{m} \sigma_k \uu_{k}$\;
      }%
      \nl \Return{$\x_N$}\;
      \caption{{SARP$(\x_0, N, m, L, [\sigma_0])$}}%
     \end{algorithm}%
     \vspace{0.255cm}
     \begin{algorithm}[H]
     %
     %
      \NoCaptionOfAlgo
      \SetAlgoVlined
      \DontPrintSemicolon
      \SetInd{0.1em}{0.5em}
      \nl $C_0 \leftarrow I_n;\, \pp_0 \leftarrow \mathbf{0}$\;
      \nl \For{$k=1$ to $N$}{
       \nl $\uu_{k} \sim \mathcal{N}(0, C_{k-1})$\;     
       \nl $(\x_k, \sigma_k) \leftarrow \mathrm{aSS}(\x_{k-1}, \uu_k, \sigma_{k-1})$\;
       \nl $\yy_{k} \leftarrow (\xx_k - \xx_{k-1})/\sigma_{k-1}$\;
       \nl \eIf{$\yy_{k} \neq \mathbf{0}$ (\textit{success})}{
       \nl $\pp_{k} \leftarrow (1-c_{\rm{c}}) \pp_{k-1} + \sqrt{c_{\rm{c}}(2-c_{\rm{c}})}\yy_{k}$ \;
       \nl $C_{k} \leftarrow (1-c_{\rm{cov}}) C_{k-1} + c_{\rm{cov}} \pp_{k} \pp_{k}^T$
       }{
       \nl $C_{k} \leftarrow C_{k-1};\, \pp_{k} \leftarrow (1-c_{\rm p})\pp_{k-1}$       
       }
      }%
      \nl \Return{$\x_N$}\;
      \caption{{(1+1)-CMA$(\x_0, N, \sigma_0, p, c_{\rm{c}}, c_{\rm{cov}})$}}%
     \end{algorithm}%
    \end{minipage}%
    }
     \hfill\null\par
    \caption{RP, EP-CMA and CMA-ES schemes.}\label{algo:all}
   \end{figure}

\textbf{SARP.} Accelerated random search schemes are fundamentally different from the simple random search schemes. Instead of generating only  one sequence of iterates, those algorithms typically maintain two or more sequences simultaneously (here essentially $\xx_k$ and $\yy_k$, see Fig.~\ref{algo:all}). Those sequences allow to store gathered knowledge on the objective function which yields better performance. In Fig.~\ref{algo:all} we present a simple version of the accelerated random search scheme proposed in~\cite{Stich:2011} and refer to it as simple accelerated random search (SARP). Like RP, SARP can (in practice) be used with exact line search oracles or with adaptive step size control, although convergence for those oracles has not been proven yet. For Nesterov's accelerated random search scheme~\cite{nesterov:11}, the expected one-step progress can be estimated as
\begin{align}
\label{eq:SARP}
\E \left[f(\x_+) \mid \x\right] \leq \left(1- (n\sqrt{\kappa})^{-1} \right) f(\x) \,,
\end{align}
where condition $\kappa = L/m$ and the two parameters $m \leq \lambda_{\rm min}(A)$ and $L \geq \lambda_{\rm max}(A)$ are required as input to the algorithm (and always provided in our numerical study). This rate is much better than~\eqref{eq:generalRP} and we hope to see that SARP attains comparable performance. SARP is not a monotone scheme, that is, the function values of the iterates are not monotonically decreasing. SARP is closely related to the first-order accelerated search scheme of Nesterov~\cite{Nesterov:2004}. This scheme also simultaneously maintains two sequences $\xx'_k$ and $\yy'_k$ of iterates (but requires access to the gradient in every iteration). For Nesterov's first-order scheme it is known~\cite[p.79]{Nesterov:2004}
that the sequence $\yy'_{k}$ obeys
\begin{align}
\begin{split}
\label{eq:ARP}
\yy'_{k+1} 
 &= \xx'_{k+1} + \beta' \left(\xx'_{k+1}- \xx'_k\right)\,,
\end{split}
\end{align}
for $\beta'=1 - 2/\sqrt{\kappa} + O(1/\kappa)$. Thus the additional $(\xx'_{k+1}-\xx'_{k})$ acts like a drift term, cf.~\cite{Polyak:1987}. For SARP with parameter $\theta'= \sqrt{1/(n^2 \kappa)}$ (only slightly different from $\theta$ in Fig.~\ref{algo:all}) the same reformulation of the update reveals
\begin{align}
\label{eq:sarp}
\yy_{k+1} = \xx_k + \beta \left(\xx_{k+1} - \xx_{k}\right)\,,
\end{align}
for $\beta= (1-\theta')/(1+\theta') = 1 - 2/(n \sqrt{\kappa}) + O(1/\kappa)$.
The main term contributing to the drift is approximately only an $1/n$-fraction of the step, accounting for the uncertainty emerging form the randomness.

\textbf{(1+1)-CMA-ES.}
In contrast to the presented Random Pursuit schemes, in CMA-ES new search points are sampled from a multivariate normal distribution $\uu_k \sim \mathcal{N}(0, C_k)$ whose parameter $C_k$ is updated in each iteration based on the evaluation of the samples. The covariance matrix can be adapted using different rank-1~\cite{Hansen:2001,igel:06} or rank-k updates~\cite{Hansen:2003}. 
In addition, the CMA-ES scheme is augmented by an auxiliary variable called evolution path that takes into account the correlation of successive means taken over a finite horizon. In~\cite{Hansen:2001,ostermeier1994}, the evolution path $\pp_k$ is updated as 
\begin{align}
\label{eq:path}
\pp_{k+1} = (1-c_{\rm c}) \pp_k + \sqrt{c_{\rm c}(1-c_{\rm c})} \uu_k\,.
\end{align}
Cumulative information about successive steps is stored in the variable $\pp_{k}$. 
We use a simplistic CMA-ES variant, closely following~\cite{Hansen:2001}, see Fig.~\ref{algo:all}. We use the simple Adaptive Step Size control aSS to determine the step size $\sigma_k$, the covariance matrix update solely uses the information of the evolution path like in~\cite{igel:06} and for simplicity we refrain from implementing any regularization features, in contrast to~\cite{igel:06}. We use the same parameters that were proposed in~\cite{igel:06} for the (1+1)-CMA-ES, namely $c_{\rm c} = 2/(n+2)$, $c_{\rm p} = 1/12$ and $c_{\rm cov}=2/(n^2+6)$.

\textbf{EP-CMA-1.}
The evolution path $\pp_k$ accumulates information over successful steps. This accumulation can be seen as a smoothing of the noisy information obtained in single steps, at the effect that the evolution path points into direction of more promising function values~\cite{Hansen:2001}. In this study, we are interested if the evolution path can be used in a similar way as the drift term in~\eqref{eq:ARP} or~\eqref{eq:sarp}, respectively, to accelerate the search. There are several ways to incorporate the evolution path $\pp_k$ into the update scheme. We suggest to use the path $\pp_k$ in the following way: in the simple random search scheme RP (equivalent to (1+1)-ES), we sample in iteration $k$ a direction from $\uu_k \sim \mathcal{N}(0, (1-c_{\rm{cov}}) I_n+c_{\rm{cov}} \pp_{k} \pp_{k}^T)$, with bias along the direction indicated by $\pp_k$. This has the effect that we only follow successful steps, but the drift imposed by the evolution path might be smaller than it ideally should be. The scheme EP-CMA-1 is detailed in Fig~\ref{algo:all}, we used $c_{\rm c}$ and $c_{\rm p}$ as above, and $c_{\rm cov} = 1/5$. This approach is similar to~\cite{sun:2013}.

\textbf{EP-CMA-$\mathbf{m}$.} The proposed EP-CMA-1 can easily be generalized to a whole family of optimization schemes by an approach presented in~\cite{Loshchilov:2014}. In EP-CMA-1, only the information stored in the current evolution path $\pp_k$ is used to bias the search direction. But we could also afford to temporarily store a small number $m$ of past $\pp_{k'}$ for $k'<k$, and use the information collectively to bias the search. As two successive evolution paths are likely highly correlated, we propose to store the evolution path only every $n^2\!/m$-th generation (and up to at most $(m-1)$ copies simultaneously). The resulting scheme is detailed in Fig.~\ref{algo:all}. We used $c_{\rm c} = 2/(n+2)$ as in CMA-ES, and for $m > 1$, $c_{\rm cov} = 2/(6+m)$ for EP-CMA-$m$. 
If implemented carefully, EP-CMA-$m$ has $\Theta(mn)$ complexity per iteration (not shown in Fig~\ref{algo:all}). For $m = n^2$, the updates of EP-CMA-$m$ are identical to the updates of (1+1)-CMA-ES, if \emph{limited} to a finite horizon of $n^2$ steps. In contrast, the low memory method proposed in~\cite{knight:2007} behaves similar to CMA-ES already for $m=n$, but has iteration complexity $\Theta(nm^2)$.

\section{Empirical Study}
\label{sec:emp}
We now present the setup of our empirical study. We focus on the following schemes: (i) the two Random Pursuit schemes with adaptive step size control (denoted as RP and SARP) and with exact line search (denoted as RP-exact and SARP-exact), (ii) the simplified (1+1)-CMA-ES and (iii) the EP-CMA-$m$ schemes as introduced in Sec.~\ref{sec:algorithms}, see Fig.~\ref{algo:all}. We use EP-CMA-$m$ with parameters $m=1,2,4,\sqrt{n},n$. This totals in 10 different schemes, all of which were implemented in MATLAB and will be made available on the authors website.

We tested the performance of all algorithms on three variants of the ellipsoidal benchmark function~\cite{Hansen:2001} and the non-convex Rosenbrock function, detailed in Table~\ref{tab:functions}. The quadratic functions were chosen in such a way that the extremal values of their spectra (1 and $L$) both agree. 
We considered the quadratic functions with parameters $L=1$\textsc{e}4 and $L=1$\textsc{e}6 each, and repeat the experiments in dimensions $n=20,40,60,80,100$.

For all experiments, initial settings were $\xx_0 = \mathbf{1}$, $\sigma_0 = 1$ and $p=0.27$ (for schemes with the aSS routine). We count the number of iterations (\#ITS) needed to decrease the function value (FVAL) below $1$\textsc{e}-9.
A graphical summary of our results can be found in Fig.~\ref{fig:all}-\ref{fig:rosen} and in the appendix.
We now proceed by discussing some of the key results.

\begin{figure}[!t]
\centering
\includegraphics[width=.49\textwidth]{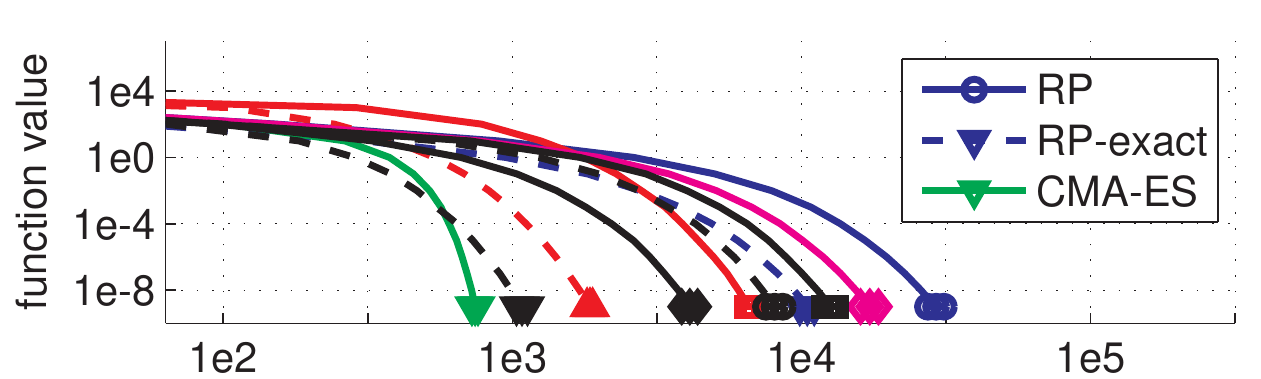}
\includegraphics[width=.49\textwidth]{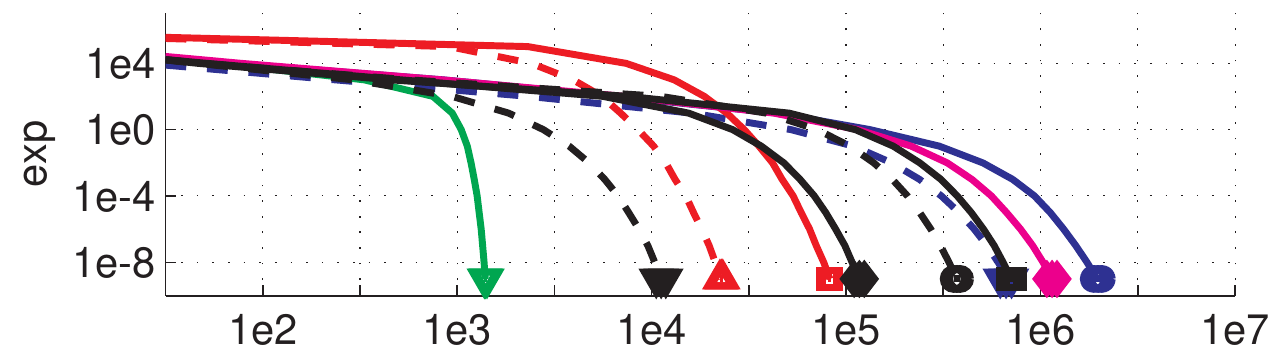}

\vspace{-1ex}
\includegraphics[width=.49\textwidth]{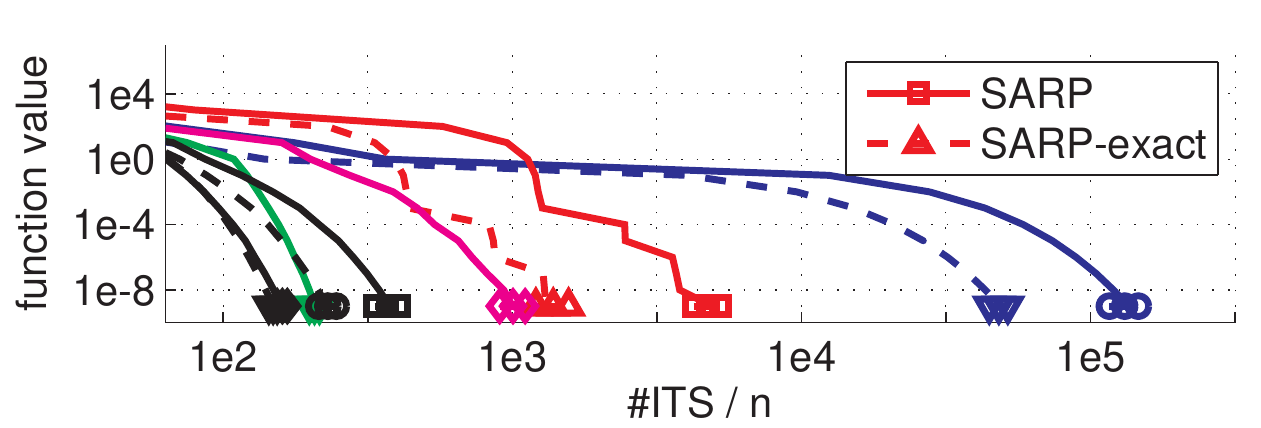}
\includegraphics[width=.49\textwidth]{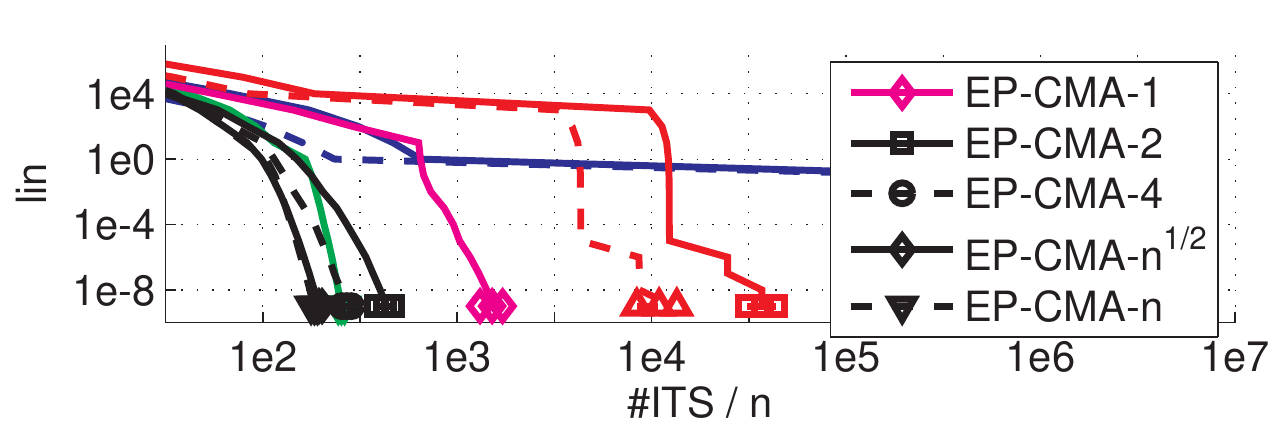}
\caption{Evolution of FVAL vs. \#ITS on $f_{\rm exp}$ (top) and $f_{\rm lin}$ (bottom) with $L=1$\textsc{e}4 (left) and $L=1$\textsc{e}6 (right) in $n=100$ dimensions. For 51 (11 for RP on $f_{\rm lin}$ with $L=1$\textsc{e}6) runs we recorded \#ITS needed to reach FVAL of 1\textsc{e}-9. The trajectory realizing the median values is depicted, mean and one standard deviation are indicated by markers. (RP on $f_{\rm lin}$ with $L=1$\textsc{e}6 reaching FVAL $<$ 1\textsc{e}-2 after 1\textsc{e}6.5$n$ \#ITS.)}
\label{fig:all}
\end{figure} 
\begin{table}[b]
  \centering
	\caption{List of benchmark functions.}
	\hrule
	\vspace{-.5em}
    \begin{align*}
    &f_{\rm exp} (\xx) = \frac{1}{2} \sum_{i=1}^n L^\frac{i-1}{n-1} x_i^2 &
    f_{\rm rosen} (\xx) = \sum_{i=1}^{n-1} \left(100 \cdot \left(x_i^2 - x_{i+1} \right)^2 + \left(x_i - 1\right)^2 \right)& \\
    &f_{\rm lin} (\xx) = \frac{1}{2} \sum_{i=1}^n \left(1 + i\frac{(L-1)}{(n-1)}\right) x_i^2     &
    f_{\rm two} (\xx) = \frac{1}{2}\sum_{i=1}^{\lfloor n/2 \rfloor} x_i^2 + \frac{L}{2} \sum_{i = \lceil n/2 \rceil} x_i^2&
    \end{align*}
    \vspace{-.5em}
 	\hrule
 \label{tab:functions}
\end{table}   
\textbf{Line search.} 
Both RP and SARP were tested with exact line search oracle and adaptive step size control. In Fig.~\ref{fig:all} we see that the exact schemes outperform their adaptive variants in $n=100$ dimensions by a factor of roughly 2-3. This pattern is observed throughout the whole benchmark in all dimensions. Thus we omit to display the results for exact line search in subsequent Figs.~\ref{fig:dimfexp}-\ref{fig:rosen}. 

\textbf{SARP vs. EP-CMA-1.} 
The picture is twofold. In Fig.~\ref{fig:all} we see 
that EP-CMA-1 outperforms SARP by a factor of roughly 5 on $f_{\rm lin}$ with $L=1$\textsc{e}4 (factor 24 for $L=1$\textsc{e}6). The smallest eigenvalue of this function is separated form the second largest by a gap of roughly $n$. Hence, knowledge of one important direction reduces the conditioning of the function by a large factor. This factor becomes smaller in higher dimension. This scaling in the dimensions is indeed observed empirically, and depicted in the appendix. 

On the other three functions, SARP performs consistently better than EP-CMA-1. On $f_{\rm exp}$ with $L=1$\textsc{e}4 in $n=100$ dimensions the factor is roughly 3, its roughly 14 for $L=1$\textsc{e}6  (Fig.~\ref{fig:all}), and exceeds 10 on both $f_{\rm two}$ and $f_{\rm rosen}$ (Fig.~\ref{fig:rosen}).
%
Considering the scaling in dimension (Figs.~\ref{fig:dimfexp}-\ref{fig:rosen}; and the appendix), %
we observe that the relative performance (\#ITS$/n$) of SARP remains constant on all four benchmark functions, as predicted by theory for a similar method~\cite{nesterov:11,lee:2013}.


\begin{figure}[!t]
\centering
\includegraphics[width=.49\textwidth]{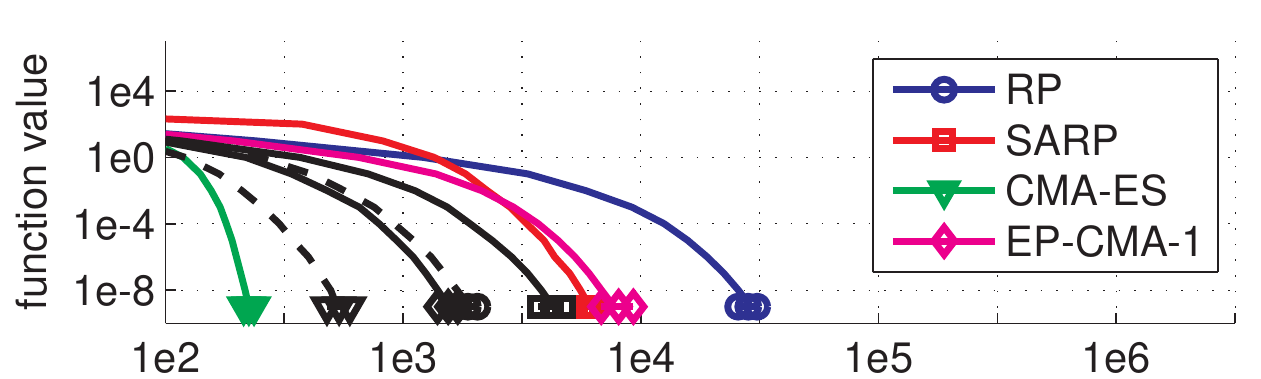}
\includegraphics[width=.49\textwidth]{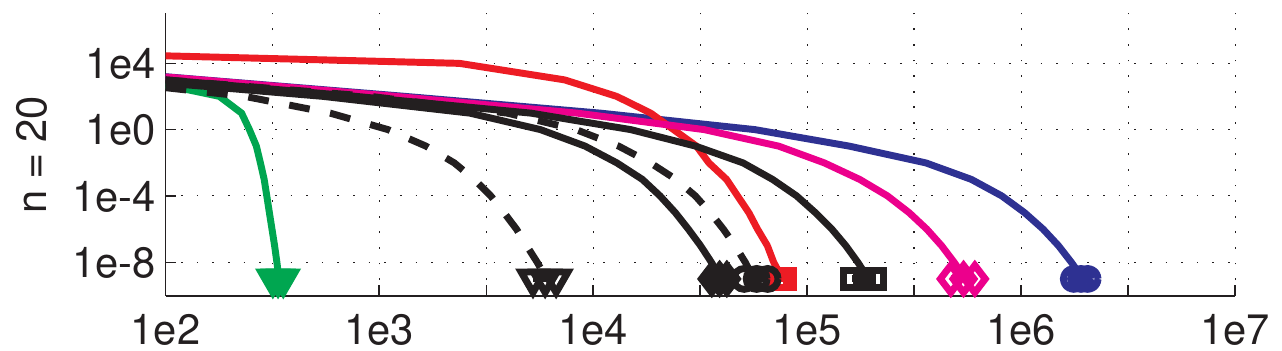}

%

\vspace{-1ex}
\includegraphics[width=.49\textwidth]{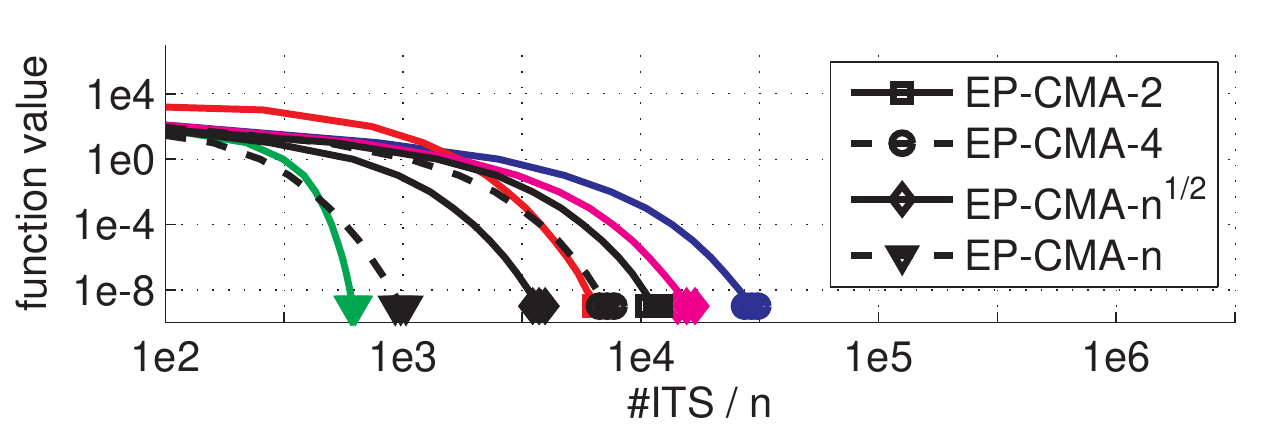}
\includegraphics[width=.49\textwidth]{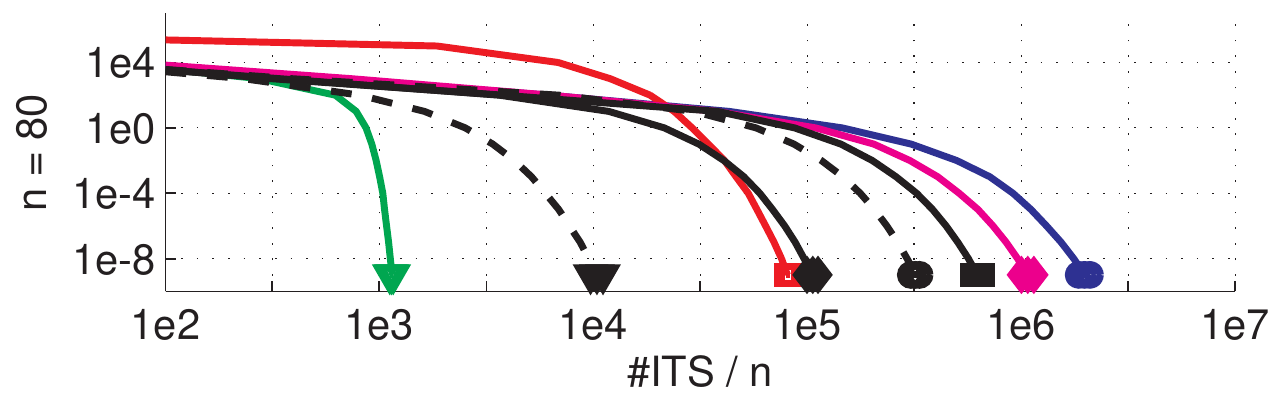}

\caption{Evolution of FVAL vs. \#ITS on $f_{\rm exp}$ with $L=1$\textsc{e}4 (left) and $L=1$\textsc{e}6 (right) in $n=20$ and $n=80$ dimensions. For 51 runs we recorded \#ITS needed to reach FVAL of 1\textsc{e}-9. The trajectory realizing the median values is depicted, mean and one standard deviation are indicated by markers.}
\label{fig:dimfexp}
\end{figure}

\textbf{EP-CMA-schemes.}
The EP-CMA-$m$ schemes consistently work better for increasing values of $m$ throughout the whole benchmark (Figs.~\ref{fig:all}-\ref{fig:rosen}). On $f_{\rm exp}$ with $L=1$\textsc{e}4 the difference in \#ITS between EP-CMA-$n$ and EP-CMA-1 is roughly a factor of 10, and 20 for $L=1$\textsc{e}6  (Fig.~\ref{fig:all}). The gap becomes gradually smaller on $f_{\rm lin}$, $f_{\rm two}$ (especially for $L=1$\textsc{e}6, see appendix), 
and is insignificant on $f_{\rm rosen}$ (Fig.~\ref{fig:rosen}).
On $f_{\rm lin}$ the EP-CMA-$m$ schemes perform extremely well, already for small $m$. EP-CMA-4 performs approximately as good as CMA-ES, for both parameters $L=1$\textsc{e}4 and $L=1$\textsc{e}6 (Fig.~\ref{fig:all}).
On $f_{\rm exp}$ in $n=100$ dimensions and parameter $L=1$\textsc{e}4, both SARP and EP-CMA-4 need about the same \#ITS. For parameter $L=1$\textsc{e}4  the performance of SARP is the same as the performance of EP-CMA\nobreakdash-$\sqrt{n}$ (Fig.~\ref{fig:all}).
On both $f_{\rm two}$ and $f_{\rm rosen}$, the EP-CMA-$m$ scheme cannot reach the performance of SARP, though on $f_{\rm rosen}$ the EP-CMA-$m$ schemes perform as good as CMA-ES (Fig.~\ref{fig:rosen}).

\textbf{CMA-ES.}
Fig.~\ref{fig:dimfexp} shows nicely the quadratic dependence of the performance of CMA-ES on the dimension $n$, see appendix 
where we report the data for all considered dimensions.
The \#ITS of the Random Pursuit schemes (RP, SARP) to reach the target accuracy increases only linearly (the relative performance (\#ITS$/n$) is constant over the dimensions).
In the dimensions $n \leq 100$ considered here, CMA-ES is the best performing scheme on $f_{\rm exp}$ (Fig.~\ref{fig:all}) and $f_{\rm two}$ (Fig.~\ref{fig:rosen}); on $f_{\rm lin}$ the EP-CMA-$m$ schemes match its performance for $m \geq 4$ (Fig.~\ref{fig:all}). A notable exception is the behavior on the non-convex $f_{\rm rosen}$, where only SARP can accelerate and the other schemes, including CMA-ES, require over 10 times more \#ITS to reach the same accuracy (Fig.~\ref{fig:rosen}).
%

\begin{figure}[!t]
\centering
\includegraphics[width=.49\textwidth]{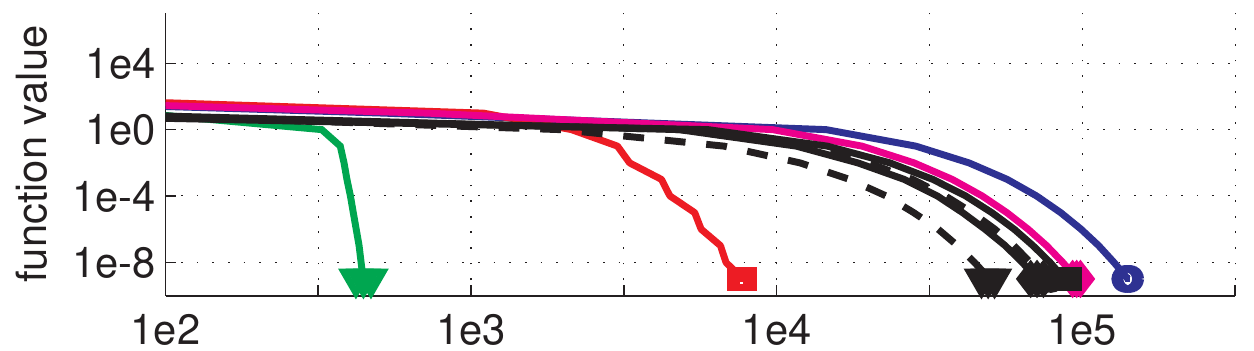}
\includegraphics[width=.49\textwidth]{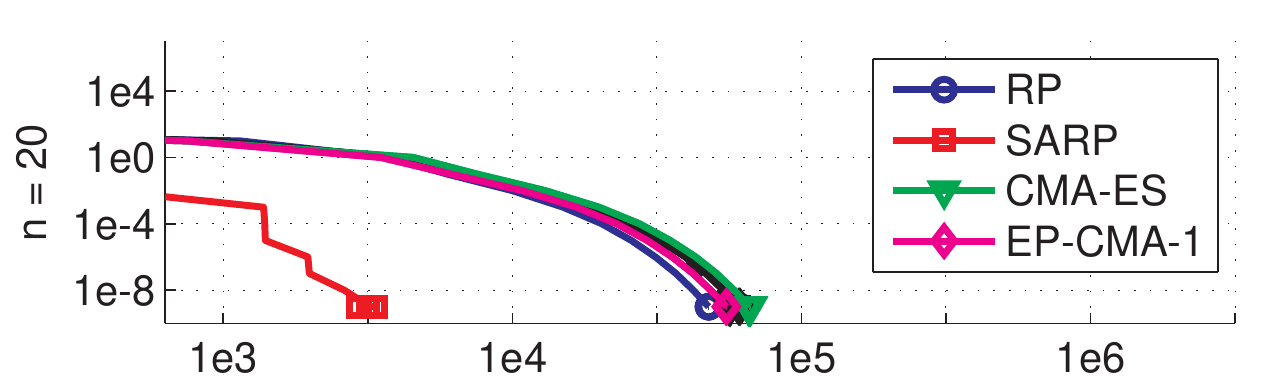}

%

\vspace{-1ex}
\includegraphics[width=.49\textwidth]{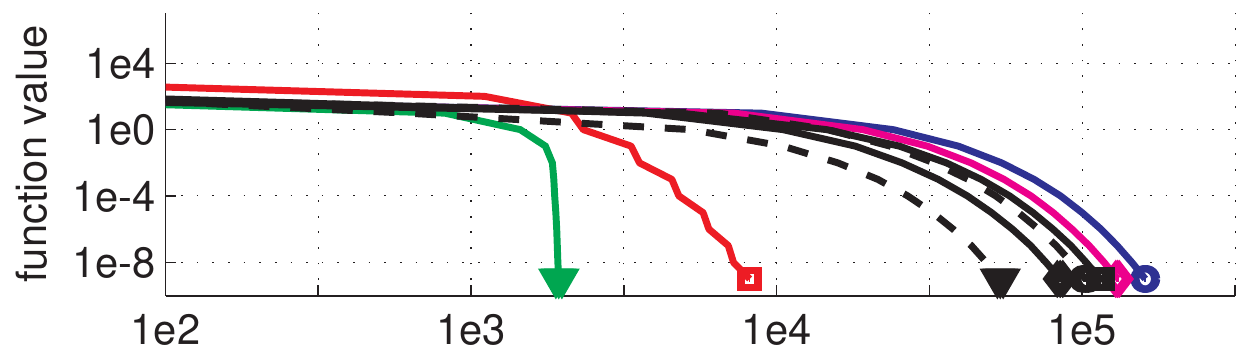}
\includegraphics[width=.49\textwidth]{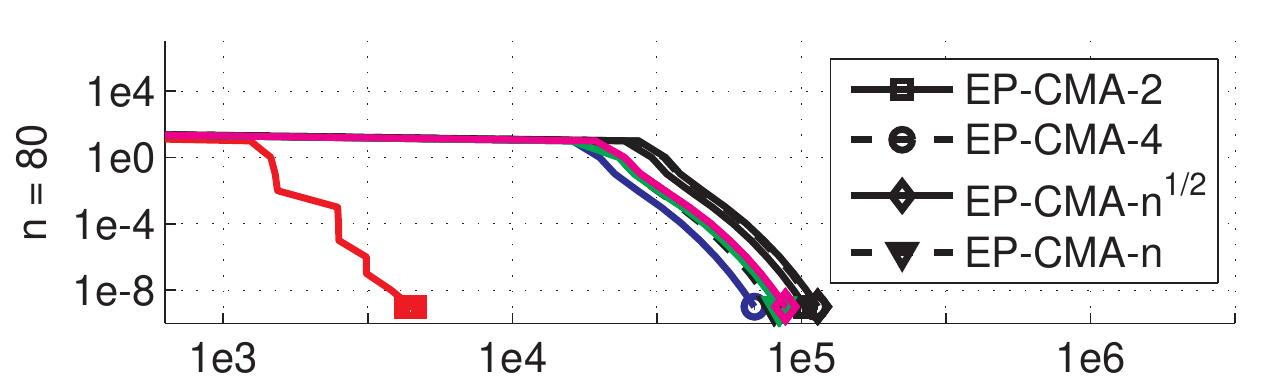}

\caption{Evolution of FVAL vs. \#ITS on $f_{\rm two}$ with $L=1$\textsc{e}4 (left) and $f_{\rm rosen}$ (right) in $n=20$ and $n=80$ dimensions. For 51 runs we recorded \#ITS needed to reach FVAL of 1\textsc{e}-9. The trajectory realizing the median values is depicted, mean and one standard deviation are indicated by markers.}
\label{fig:rosen}
\end{figure}
\section{Discussion and Conclusions}
\label{sec:disc}
In this contribution we emphasize the importance of accelerated random gradient schemes~\cite{nesterov:11,lee:2013}.
Each iteration in SARP has only linear complexity, yet the scheme takes correlations between successive iterates into account.
In CMA-ES, such correlations are collected in the evolution path~\cite{Hansen:2001,ostermeier1994} and stored in the covariance matrix. This requires $\Theta(n^2)$ simple operations per iteration. 
The proposed EP-CMA-1 uses as well the information of the evolution path to bias the search, but does not store a full-rank covariance matrix. 

\textbf{Line search.}
We empirically tested two Random Pursuit algorithms with an exact line search oracle.
Such an oracle is in general not available for general black-box optimization problems and the line search must for instance be implemented as bisection search~(cf.~\cite{Stich:2011,Stich:2012a}) at the expense of additional function evaluations per iteration. The empirical data shows that both Random Pursuit algorithms do perform well if a simple adaptive step size scheme is used instead of the line search. This makes both schemes (especially SARP) promising candidates for black-box optimization, also in high dimension, as the runtime scales only linearly with the dimension. 
Up to our knowledge, no experimental results for SARP with adaptive step size have been published yet (the authors in~\cite{Stich:2011} considered a line search with high accuracy, almost like SARP-exact). 

\textbf{Acceleration in EP schemes.} 
Our empirical results show that the sole use of the evolution path
can lead to astonishingly good performance---depending on the problem and its eigenvalue spectrum.
The speed-up of EP-CMA-1 on $f_{\rm lin}$ can be explained by the fact that the condition number of the problem drops once 
the algorithm has learned the most insensitive direction. 
Hence, the acceleration can be explained by formula~\eqref{eq:generalRP} rather than~\eqref{eq:SARP}. 
For SARP the situation is more promising. 
The data indicates that the convergence on $f_{\rm exp}$ and $f_{\rm two}$ is as described in~\eqref{eq:SARP}. The same seems to be true on the non-convex $f_{\rm rosen}$ where SARP needs an order of magnitude less \#ITS than all other schemes, including CMA-ES. Only on
$f_{\rm lin}$ this does not to hold, as SARP is only one order of magnitude faster than RP.  
Consider the update~\eqref{eq:sarp}. By expansion we obtain
\begin{align}
\yy_{k+1} 
 &= \xx_k + \beta \left(\sigma_{k+1} \uu_{k+1} + \beta \left(\xx_{k} - \xx_{k-1}\right) \right) = \xx_k + \sum_{i=1}^{k+1} \beta^{k+2-i} \sigma_i \uu_i\,.
\end{align}
We see that the drift is a weighted average of the previous steps $\sigma_i \uu_i$. The discount factor $\beta$ is the expected convergence rate. Therefore, the influence of a step $\sigma_i \uu_i$ on $\yy_{k+1}$ is roughly the same for all $i=1,\dots, k$. In contrast to this, the evolution path $\pp_k$ stores only information of the directions of the last steps (but no step sizes). The discount factor is approximately $1-2/n$. Although the evolution path $\pp_k$ is a cumulation of all old steps, the weigh of old steps is exponentially small compared to the influence of the newest steps. We might conclude that the mechanism of  accelerated random schemes like SARP is therefore inherently different to the concept of the evolution path, supporting reports in~\cite{Stich:2012a}. However, we cannot rule out the possibility, that with a different choice of internal parameters of EP-CMA-1 the difference to SARP could be reduced.

\textbf{Limited Memory schemes.}
The performance of the  proposed EP-CMA-$m$ schemes 
uniformly increases for larger parameters $m$, as well as the complexity of each single iteration. 
An optimal trade-off for the parameter $m$ has to be found, depending on the dimension $n$ and the cost of individual function evaluations. 
The data shows that the EP-CMA-$m$ schemes can dramatically improve the performance of simple random search already for small values of $m$. The speed-up depends crucially on the eigenvalue spectra of the objective function. 
It seems that these schemes can not reach the performance of the related variants in~\cite{Loshchilov:2014}.


We generally conclude, that the here
proposed algorithmic schemes with linear iteration complexity could
be a promising way to handle high dimensional black-box optimization problems.
However, the empirical data suggest that there is an intrinsic limitation for the EP schemes, as they depend on the eigenvalue spectrum of the objective function. 
%
This behavior is not observed for SARP. We like to advocate that features of accelerated schemes (like SARP) should therefore be taken seriously into account when facing high dimensional problems.


%
\bibliographystyle{splncs}
\bibliography{ppsn2014}

\newpage
\appendix
\section{Appendix}
~
\vfill
\begin{figure}[h!]
\centering
\includegraphics[width=.49\textwidth]{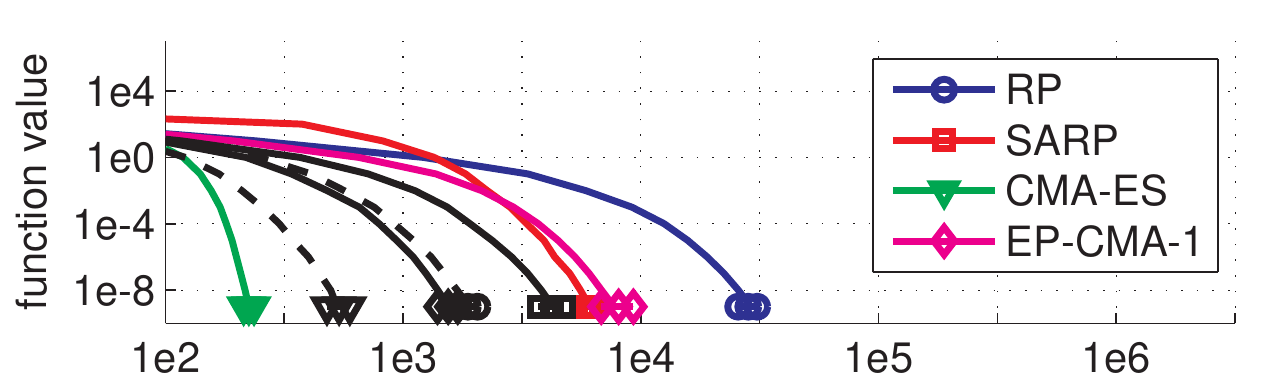}
\includegraphics[width=.49\textwidth]{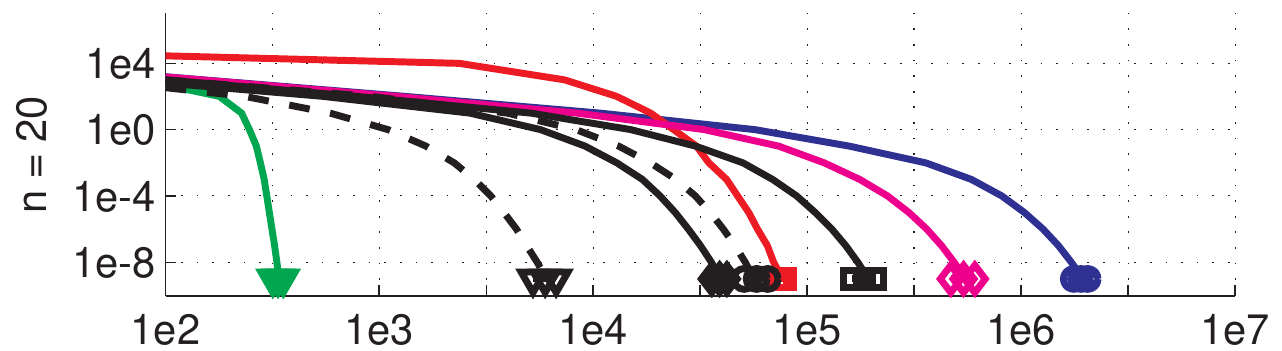}

\vspace{-1ex}
\includegraphics[width=.49\textwidth]{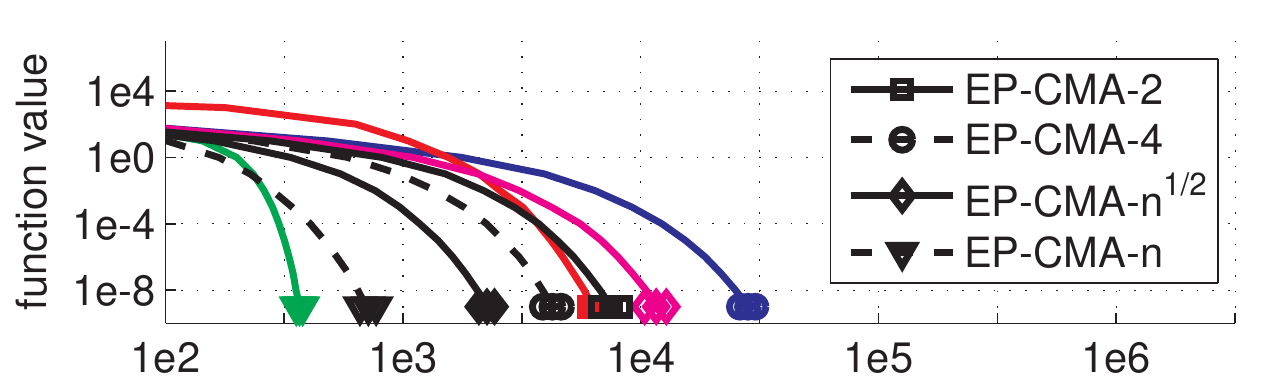}
\includegraphics[width=.49\textwidth]{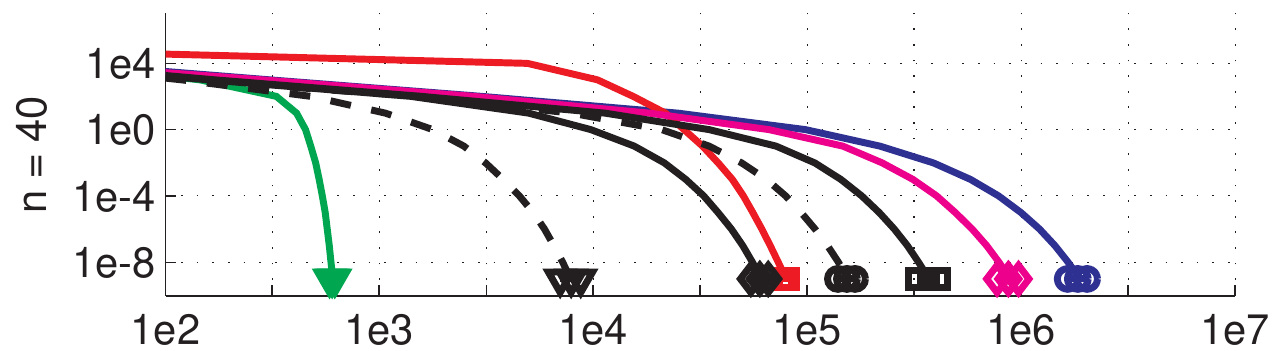}

\vspace{-1ex}
\includegraphics[width=.49\textwidth]{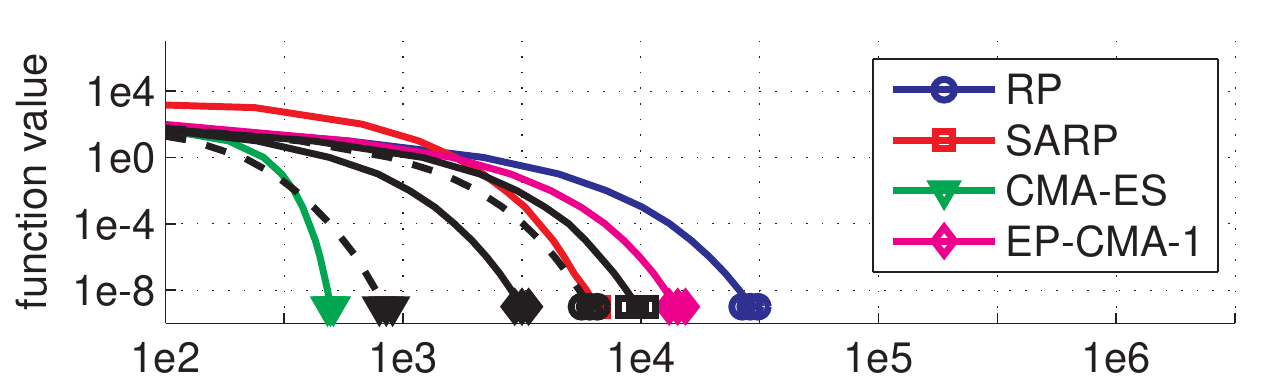}
\includegraphics[width=.49\textwidth]{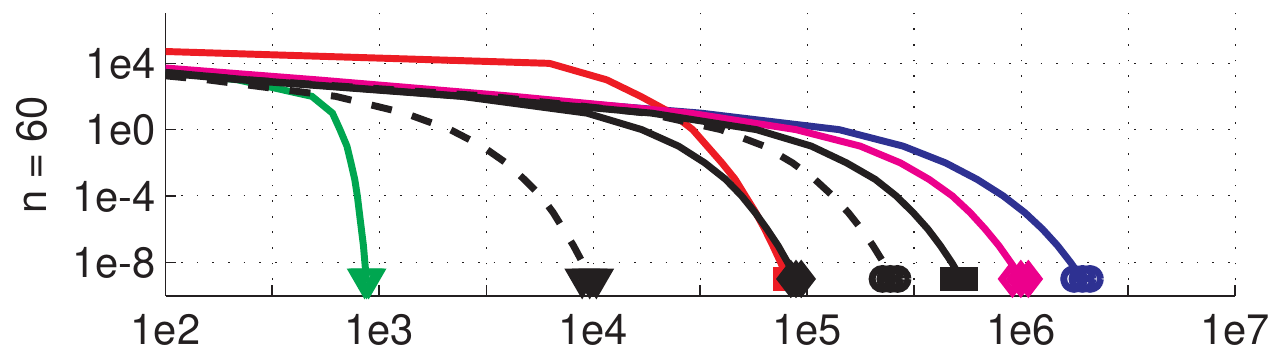}

\vspace{-1ex}
\includegraphics[width=.49\textwidth]{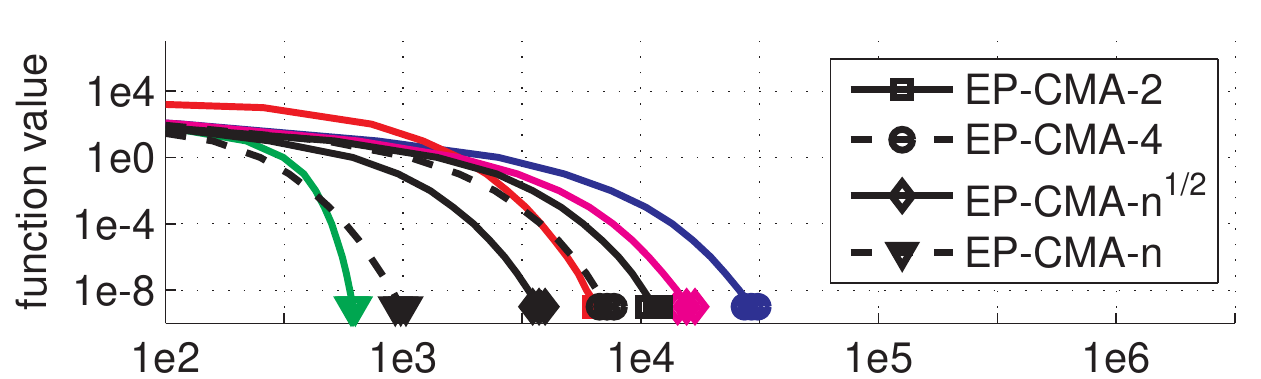}
\includegraphics[width=.49\textwidth]{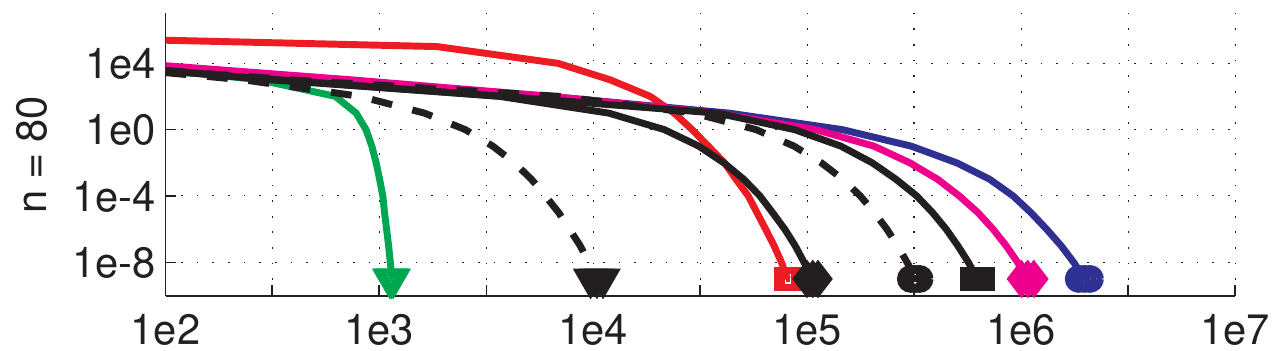}

\vspace{-1ex}
\includegraphics[width=.49\textwidth]{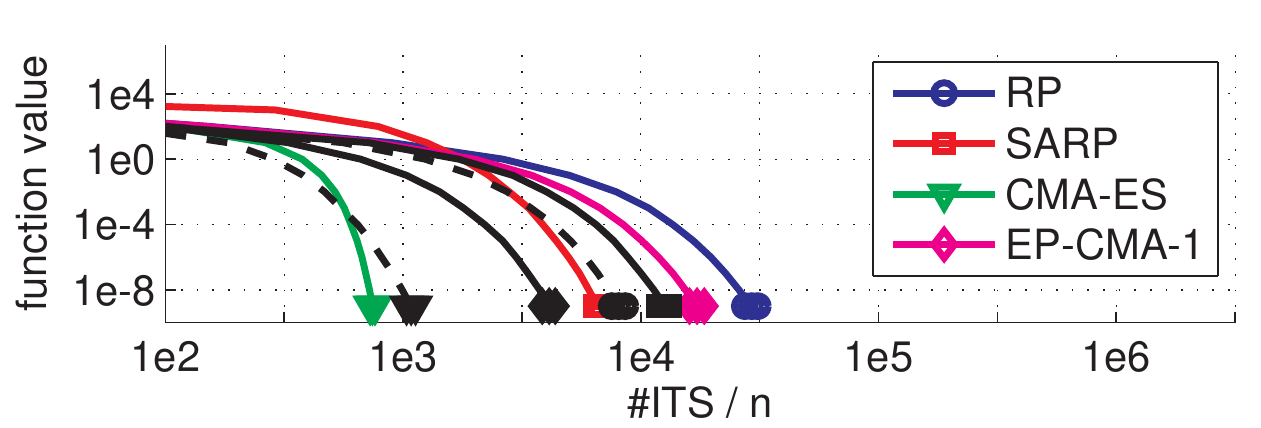}
\includegraphics[width=.49\textwidth]{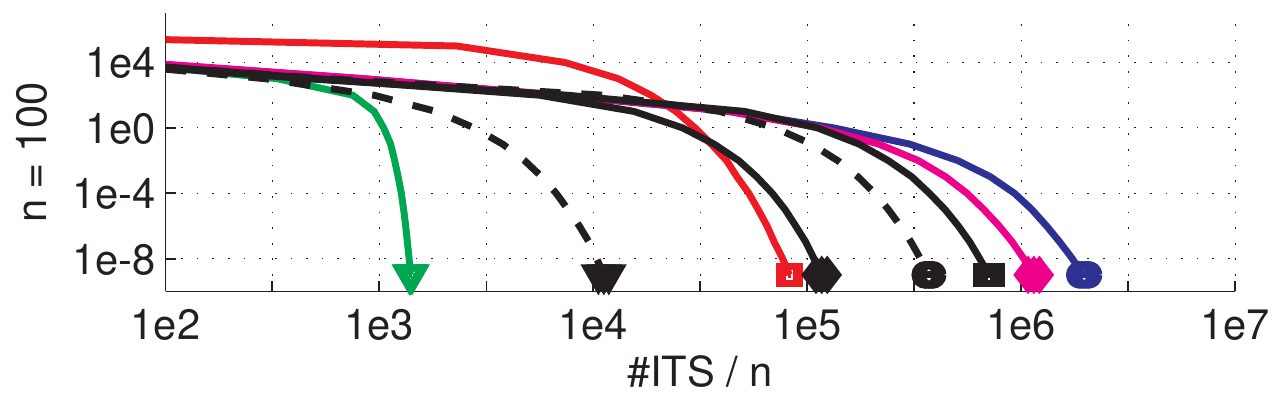}

\caption{Evolution of FVAL vs. \#ITS on $f_{\rm exp}$ with $L=1$\textsc{e}4 (left) and $L=1$\textsc{e}6 (right) in $n=\{20,40,60,80,100\}$ dimensions. For 51 runs we recorded \#ITS needed to reach FVAL of 1\textsc{e}-9. The trajectory realizing the median values is depicted, mean and one standard deviation are indicated by markers.}
\label{fig:supplement:exp}
\end{figure}
\vfill

\begin{figure}[h!]
\centering
\includegraphics[width=.49\textwidth]{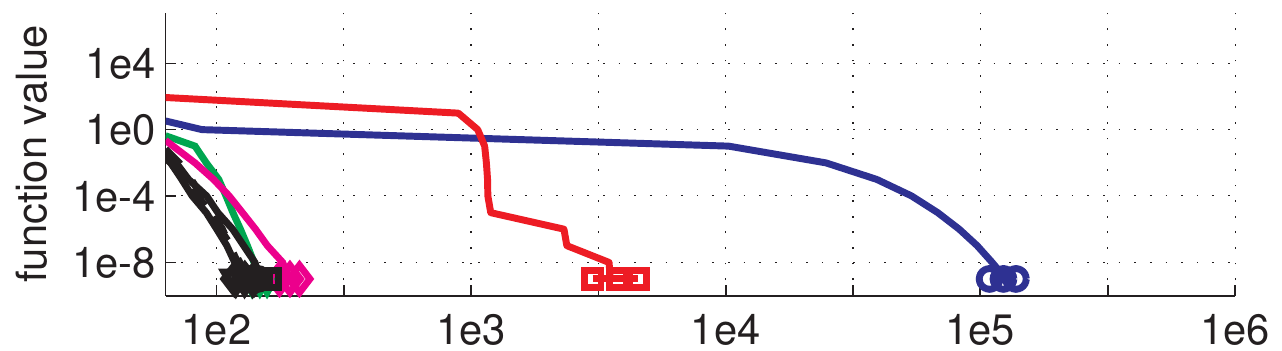}
\includegraphics[width=.49\textwidth]{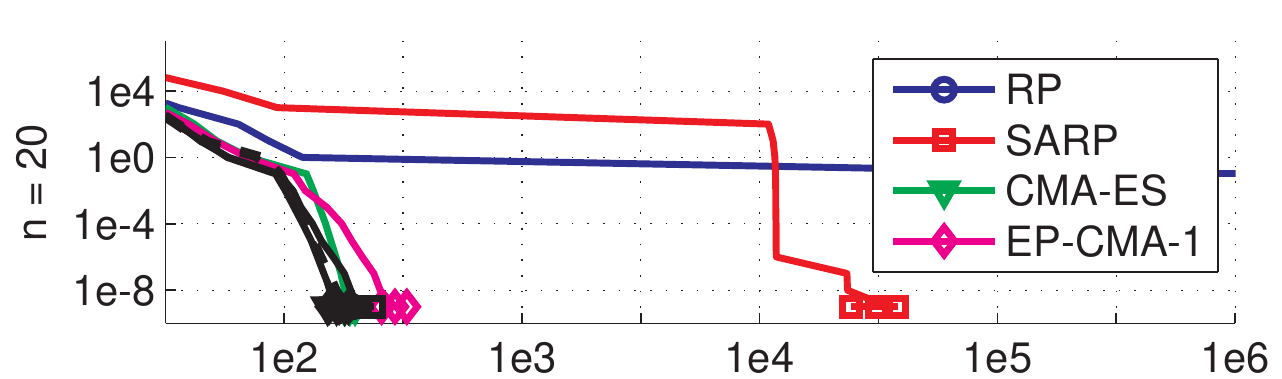}

\includegraphics[width=.49\textwidth]{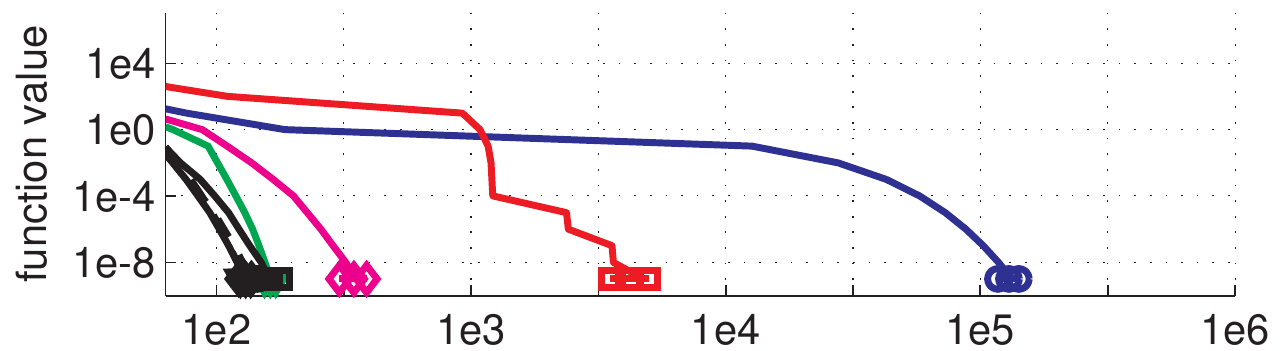}
\includegraphics[width=.49\textwidth]{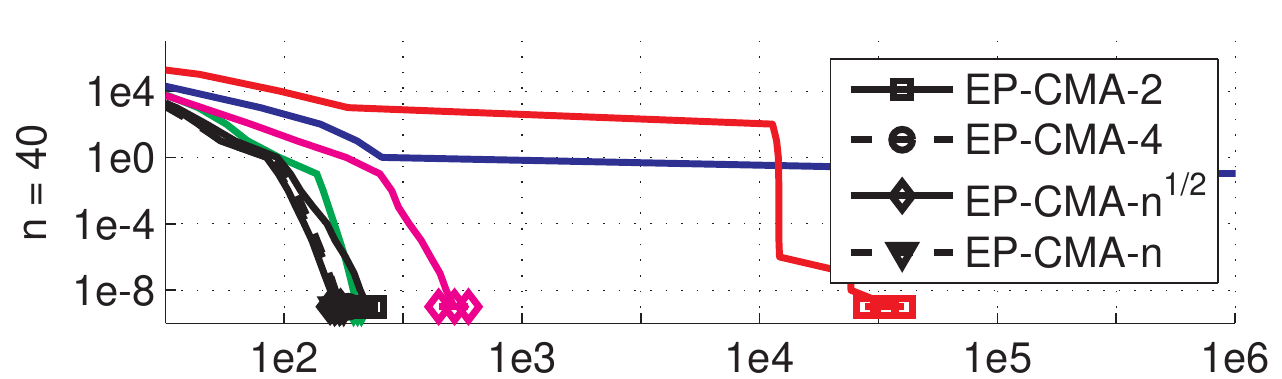}

\includegraphics[width=.49\textwidth]{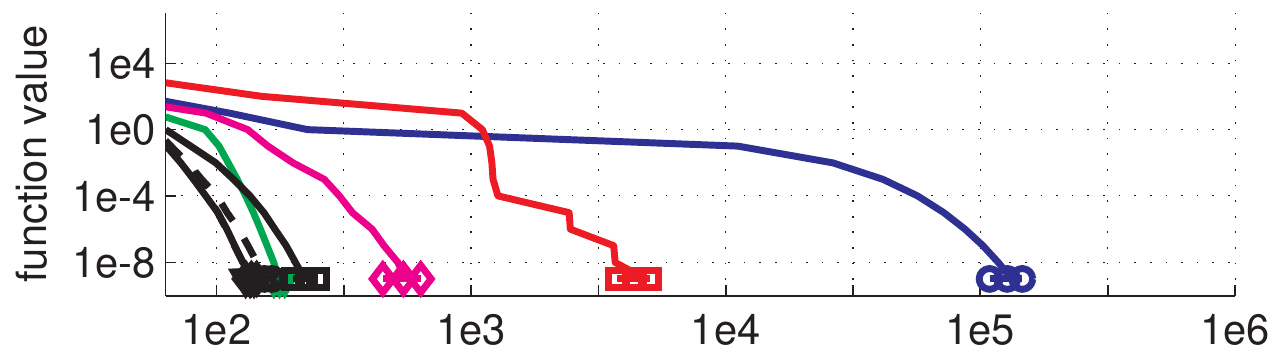}
\includegraphics[width=.49\textwidth]{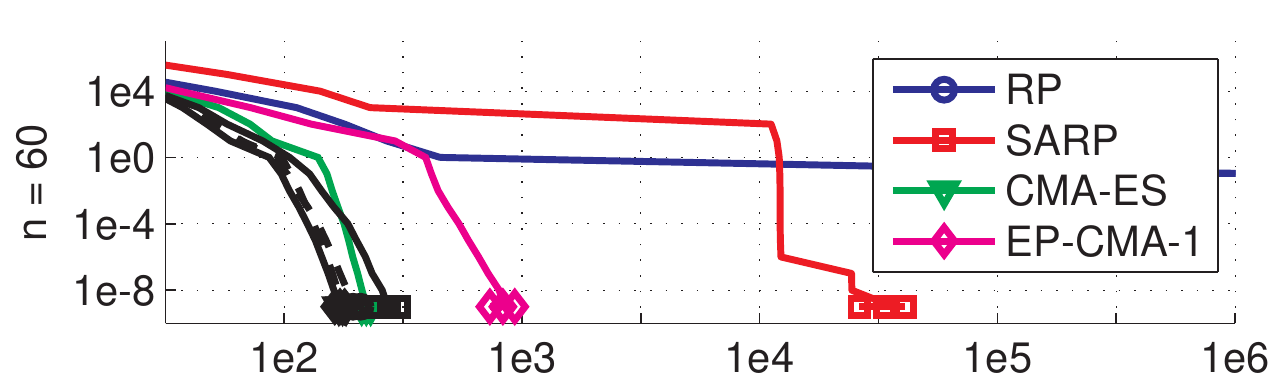}

\includegraphics[width=.49\textwidth]{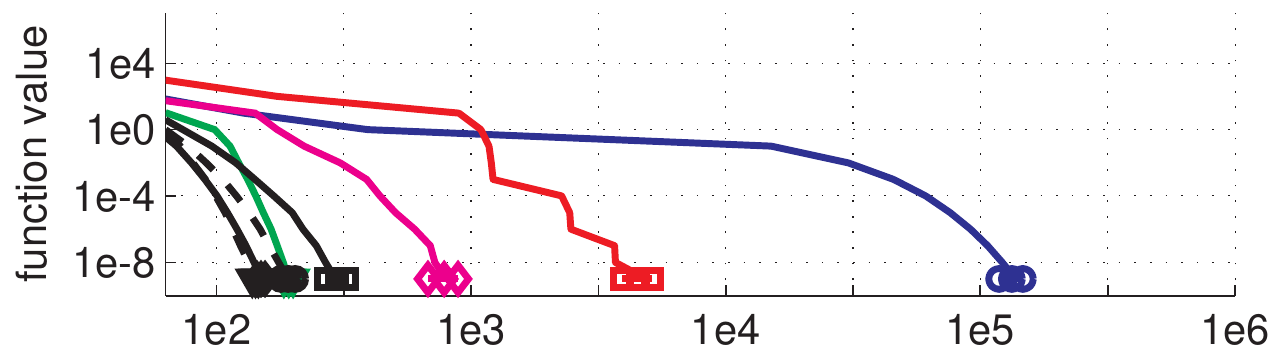}
\includegraphics[width=.49\textwidth]{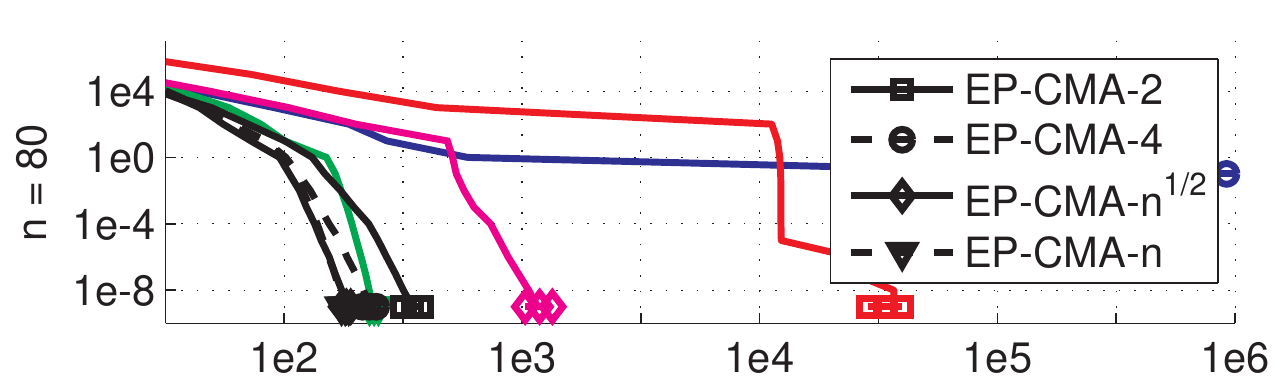}

\includegraphics[width=.49\textwidth]{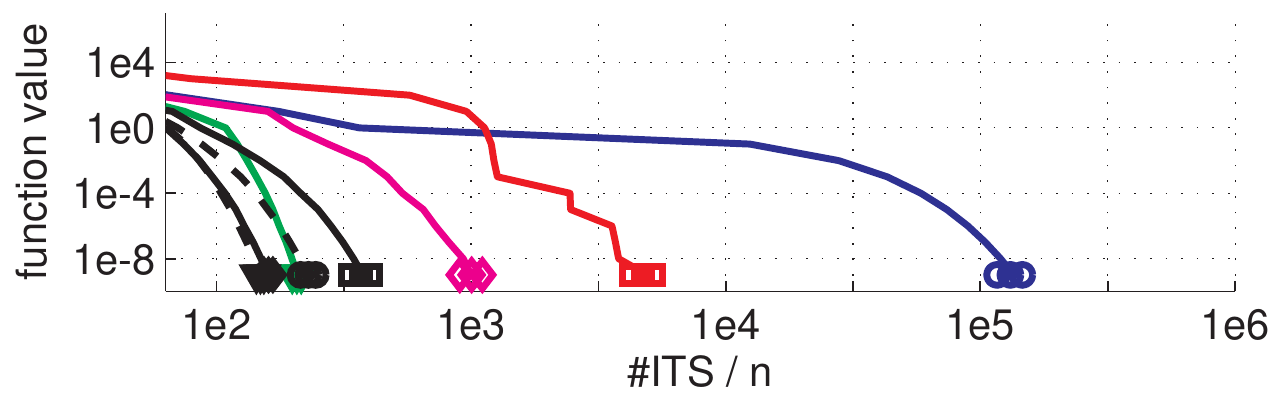}
\includegraphics[width=.49\textwidth]{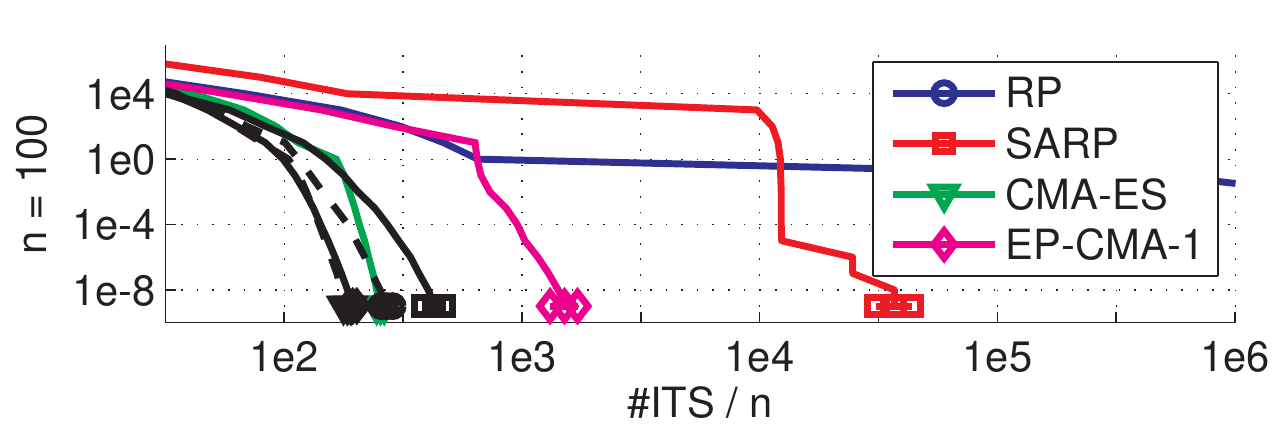}

\caption{Evolution of FVAL vs. \#ITS on $f_{\rm lin}$ with $L=1$\textsc{e}4 (left) and $L=1$\textsc{e}6 (right) in $n=\{20,40,60,80,100\}$ dimensions. For 51 runs (11 for RP and $L=1$\textsc{e}6) we recorded \#ITS needed to reach FVAL of 1\textsc{e}-9. The trajectory realizing the median values is depicted, mean and one standard deviation are indicated by markers.}
\label{fig:supplement:lin}
\end{figure} 

\begin{figure}[h!]
\centering
\includegraphics[width=.49\textwidth]{matlab/plots_supplement/ftwo4.dim_20.all=1.ls=0-eps-converted-to.pdf}
\includegraphics[width=.49\textwidth]{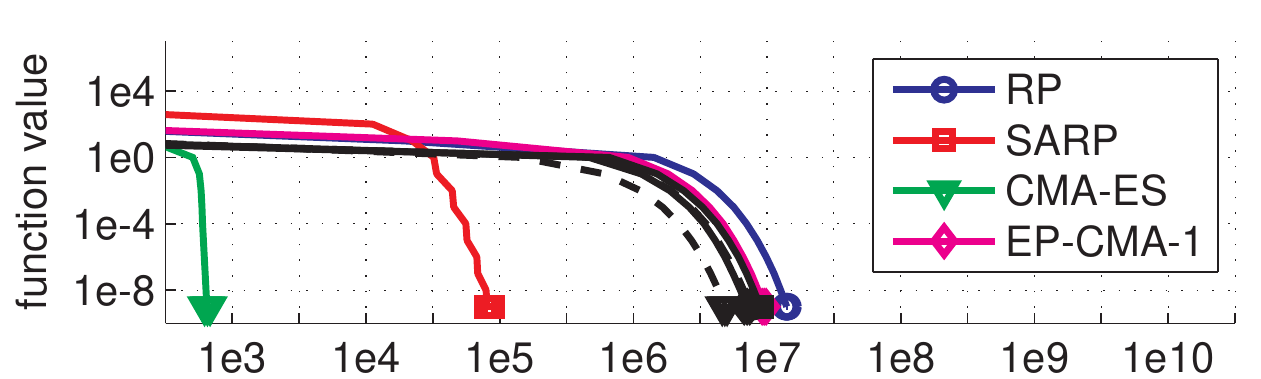}

\includegraphics[width=.49\textwidth]{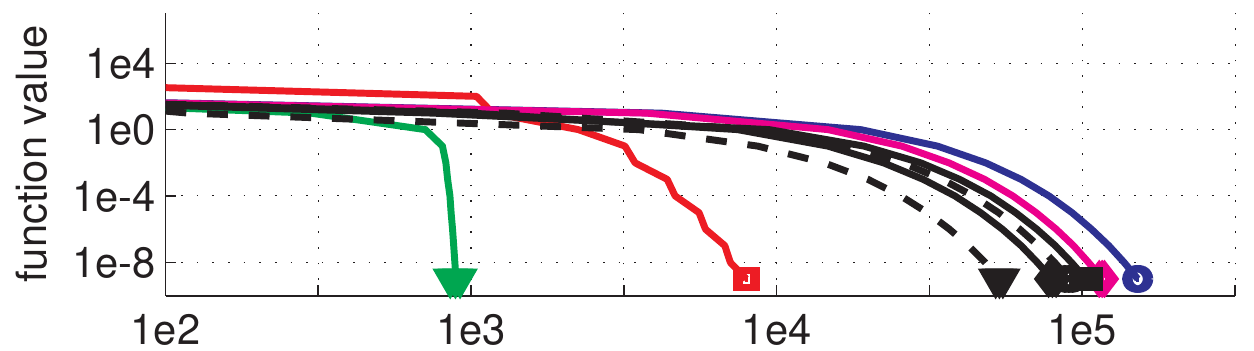}
\includegraphics[width=.49\textwidth]{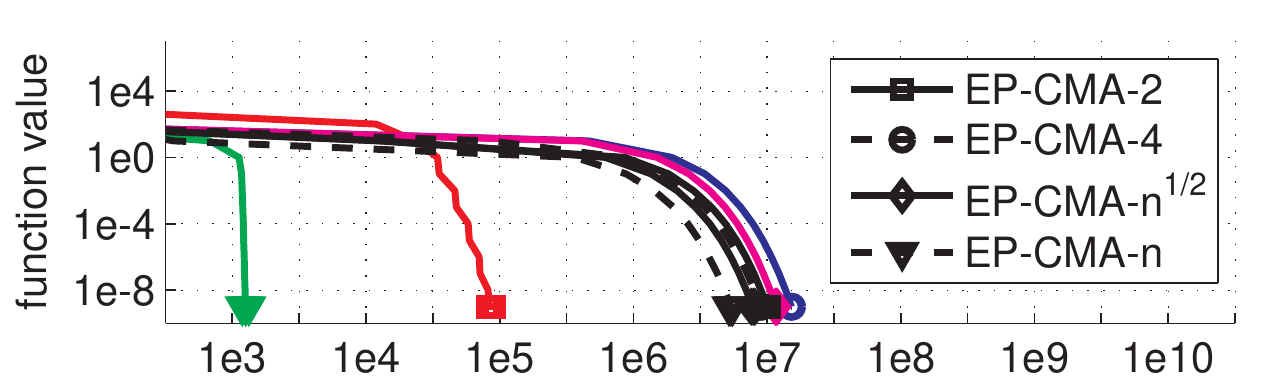}

\includegraphics[width=.49\textwidth]{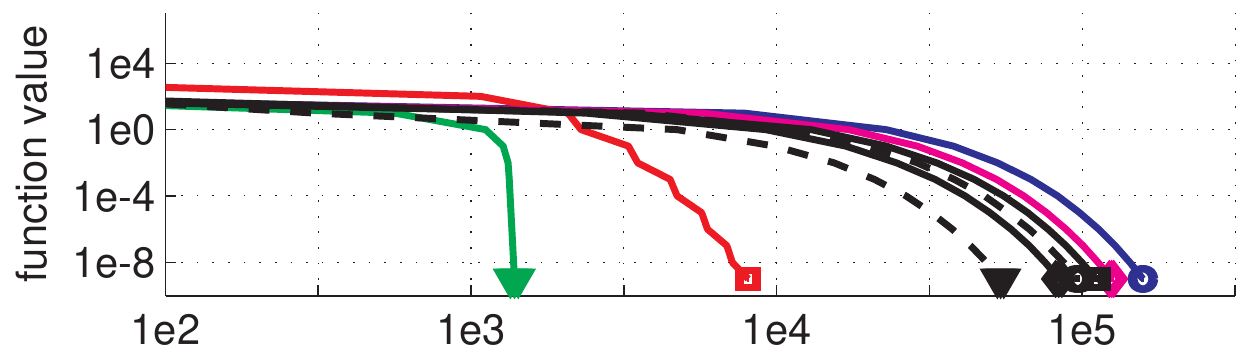}
\includegraphics[width=.49\textwidth]{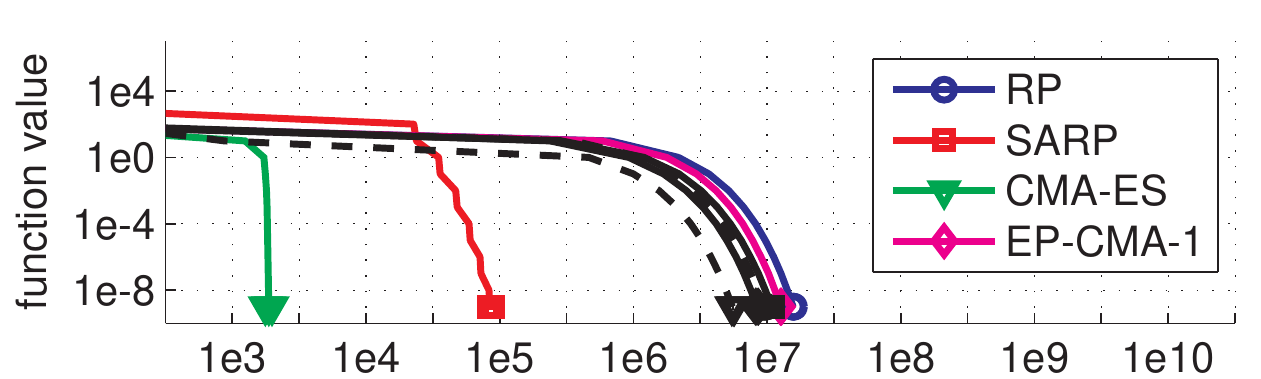}

\includegraphics[width=.49\textwidth]{matlab/plots_supplement/ftwo4.dim_80.all=1.ls=0-eps-converted-to.pdf}
\includegraphics[width=.49\textwidth]{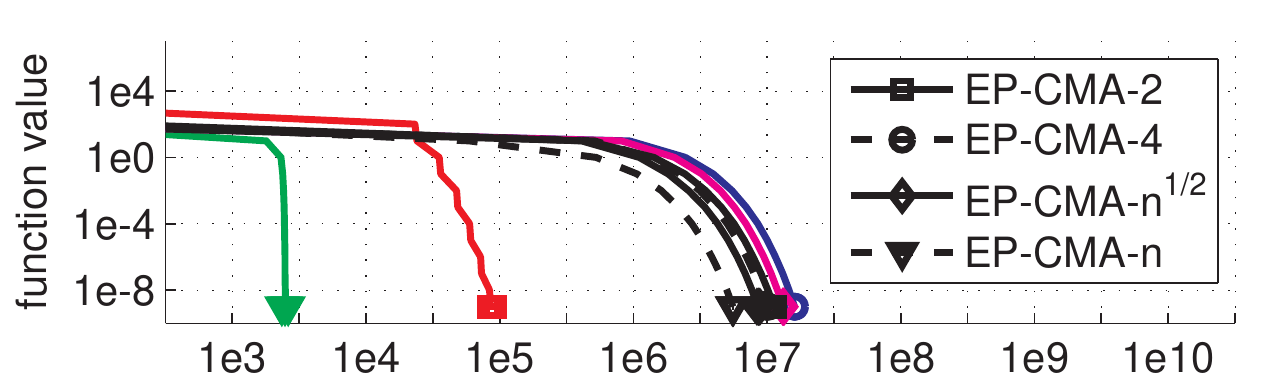}

\includegraphics[width=.49\textwidth]{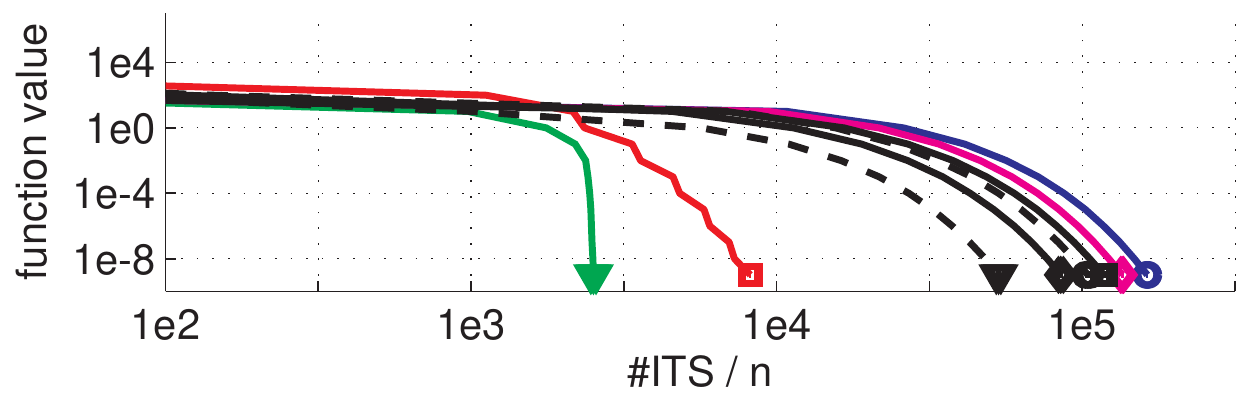}
\includegraphics[width=.49\textwidth]{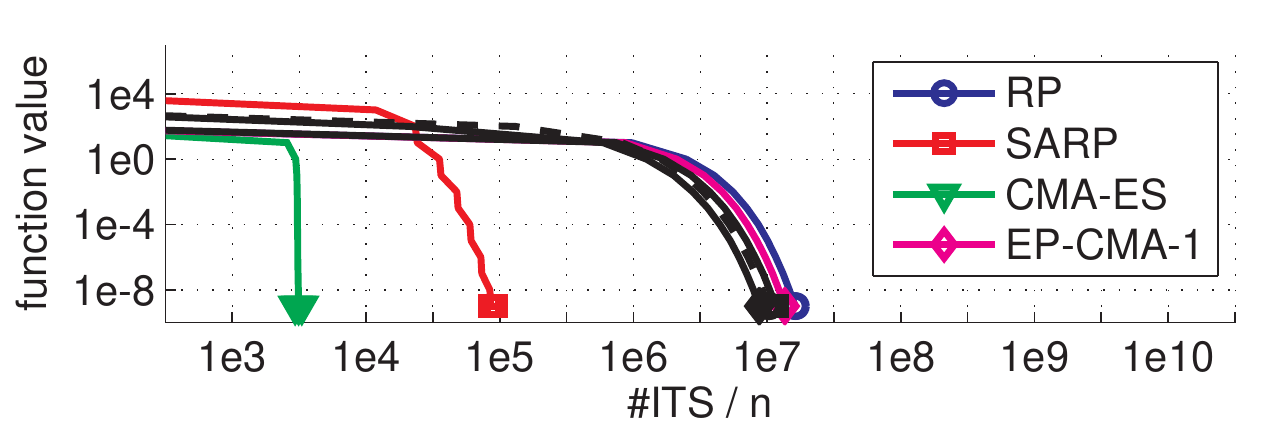}

\caption{Evolution of FVAL vs. \#ITS on $f_{\rm two}$ with $L=1$\textsc{e}4 (left) and $L=1$\textsc{e}6 (right) in $n=\{20,40,60,80,100\}$ dimensions. For 51 runs (see below) we recorded \#ITS needed to reach FVAL of 1\textsc{e}-9. The trajectory realizing the median values is depicted, mean and one standard deviation are indicated by markers.
(For $L=1$\textsc{e}6 in dimension $n=100$: only 11 runs of all schemes, except EP-CMA-100 which is omitted.)}
\label{fig:supplement:two}
\end{figure} 

\begin{figure}[h!]
\centering
\includegraphics[width=.49\textwidth]{matlab/plots_supplement/frosen.dim_20.all=1.ls=0-eps-converted-to.pdf}

\includegraphics[width=.49\textwidth]{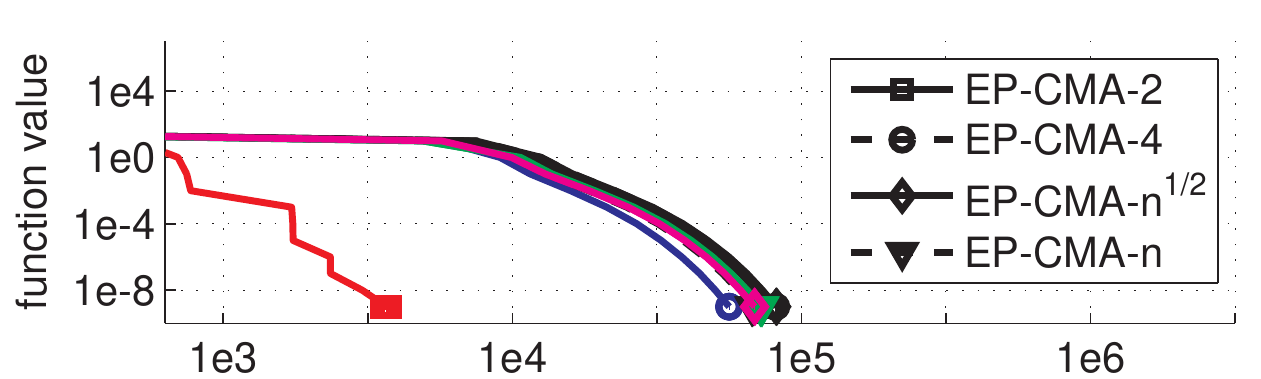}

\includegraphics[width=.49\textwidth]{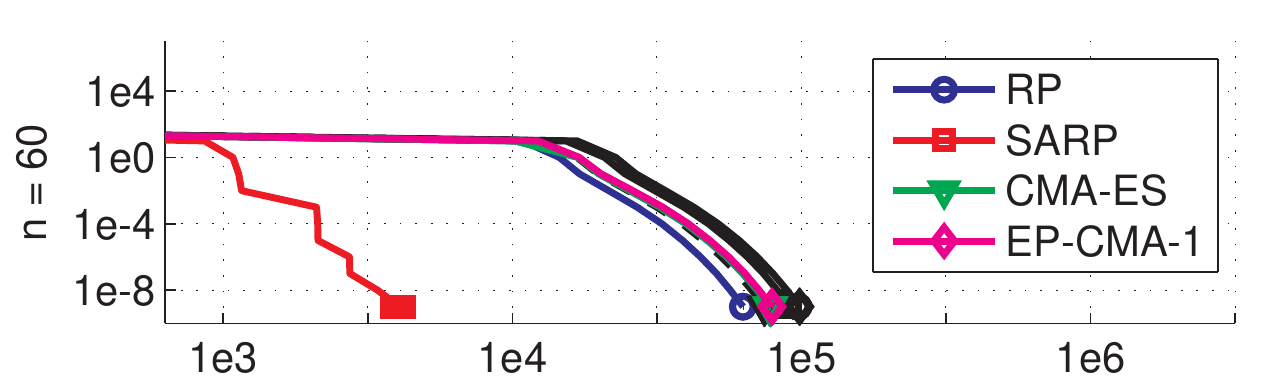}

\includegraphics[width=.49\textwidth]{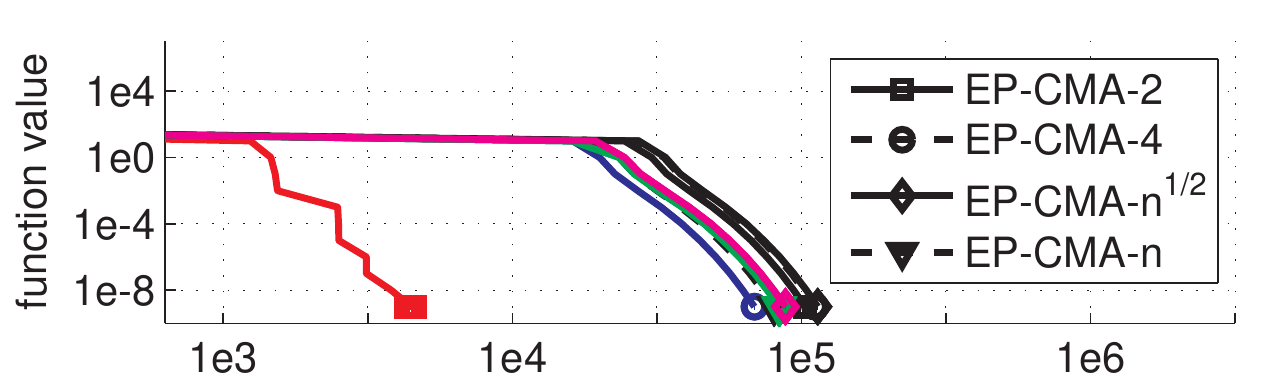}

\includegraphics[width=.49\textwidth]{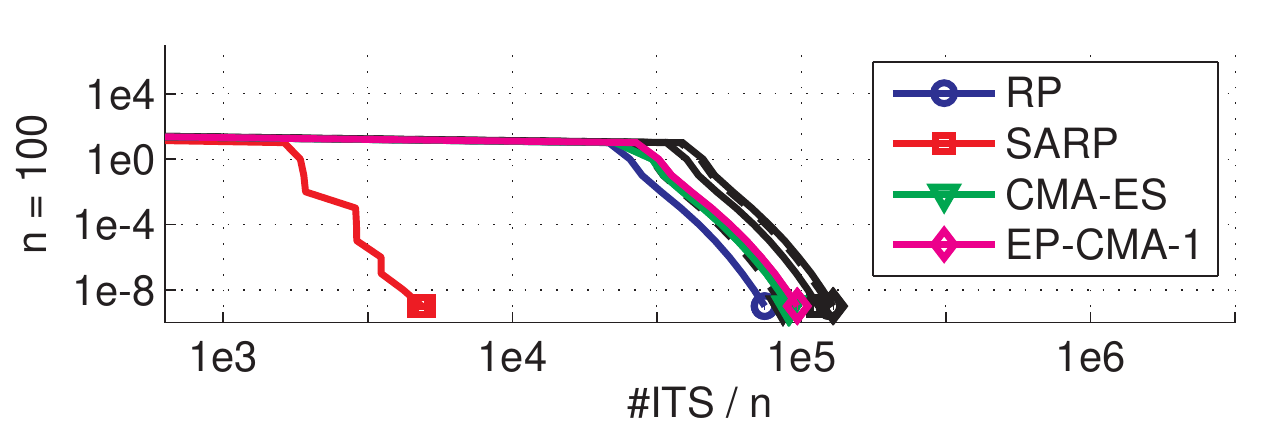}

\caption{Evolution of FVAL vs. \#ITS on $f_{\rm rosen}$ in $n=\{20,40,60,80,100\}$ dimensions. For 51 runs we recorded \#ITS needed to reach FVAL of 1\textsc{e}-9. The trajectory realizing the median values is depicted, mean and one standard deviation are indicated by markers.}
\label{fig:supplement:rosen}
\end{figure}

\end{document}